\documentclass{amsart}
\usepackage{amssymb,latexsym,graphicx}

\DeclareMathOperator{\tr}{tr}
\DeclareMathOperator{\id}{id}
\DeclareMathOperator{\supp}{supp}
\DeclareMathOperator{\im}{im}
\DeclareMathOperator{\coker}{coker}

\def\L{\mathcal{L}}
\def\J{\mathcal{J}}
\def\P{\mathcal{P}}

\def\O{\mathcal{O}}
\def\F{{\mathcal{F}}}
\def\TT{{\mathcal{T}}}
\def\injlim{\lim\limits_\rightarrow}

\def\T{\mathbb T}
\def\R{\mathbb R}
\def\Z{\mathbb Z}
\def\N{\mathbb N}
\def\C{\mathbb C}
\def\Q{\mathbb Q}

\def\zeroOne{zero-one}

\def\htau{\widehat \tau}
\def\pwm{piecewise monotonic}

\def\mapright#1{\smash{\mathop{\longrightarrow}\limits^{#1}}}
\def\mapdown#1{\downarrow\rlap{$\vcenter{\hbox{$\scriptstyle#1$}}$}}
\def\mapup#1{\uparrow\rlap{$\vcenter{\hbox{$\scriptstyle#1$}}$}}
\def\mapleft#1{\smash{\mathop{\longleftarrow}\limits^{#1}}}

\theoremstyle{plain}

\newtheorem{theorem}{Theorem}[section]
\newtheorem{proposition}[theorem]{Proposition}
\newtheorem{lemma}[theorem]{Lemma}
\newtheorem{corollary}[theorem]{Corollary}

\newtheorem*{IntroTheorem}{Theorem}
\def\Theo[#1]#2{\begin{IntroTheorem}[#1]#2\end{IntroTheorem}}

\theoremstyle{definition}
\newtheorem{definition}[theorem]{Definition}
\newtheorem{example}[theorem]{Example}

\theoremstyle{remark}
\newtheorem*{notation}{Notation}

\newif\ifproofing
\proofingfalse
\newcommand{\note}[1]{\ifproofing\marginpar{#1}\fi}

\newcommand{\prop}[2]{\begin{proposition}\label{#1}\note{#1}#2\end{proposition}}
\newcommand{\lem}[2]{\begin{lemma}\label{#1}\note{#1}#2\end{lemma}}
\newcommand{\cor}[2]{\begin{corollary}\label{#1}\note{#1}#2\end{corollary}}
\newcommand{\theo}[2]{\begin{theorem}\label{#1}\note{#1}#2\end{theorem}}
\newcommand{\defi}[2]{\begin{definition}\label{#1}\note{#1}#2\end{definition}}
\newcommand{\exe}[2]{\begin{example}\label{#1}\note{#1}#2\end{example}}
\newcommand{\notate}[1]{\begin{notation}#1\end{notation}}

\newcommand{\prooff}[1]{\begin{proof}#1\end{proof}}

\begin{document}

\title{C*-algebras associated with interval maps}
\author{Valentin Deaconu}
\address{University of Nevada\\Reno, Nevada}
\email{vdeaconu@unr.edu}
\author{Fred Shultz}
\address{Wellesley College\\Wellesley, Massachusetts 02481}
\email{fshultz@wellesley.edu}
\subjclass[2000]{Primary 46L80; Secondary 37E05}
\keywords{dimension group, interval map, piecewise monotonic, unimodal map,
tent map, Markov map, $\beta$-shift, interval exchange map, C*-algebra, Cuntz-Krieger algebra, Matsumoto algebra}

\begin{abstract}
For each piecewise monotonic map $\tau$ of $[0,1]$, we associate a
pair of C*-algebras $F_\tau$ and $O_\tau$ and calculate their
K-groups. The algebra $F_\tau$ is an AI-algebra.  We characterize
when $F_\tau$ and $O_\tau$ are simple.  In those cases, $F_\tau$
has a unique trace, and $O_\tau$ is purely infinite with a unique
KMS-state. In the case that $\tau$ is Markov, these algebras
include the Cuntz-Krieger algebras $O_A$, and the associated
AF-algebras $F_A$. Other examples for which the K-groups are
computed include tent maps, quadratic maps, multimodal maps,
interval exchange maps, and $\beta$-transformations.  For the case
of interval exchange maps and of $\beta$-transformations, the
C*-algebra $O_\tau$ coincides with  the algebras defined by Putnam
and Katayama-Matsumoto-Watatani respectively.
\end{abstract}

\maketitle

\section{Introduction}

\subsection*{Motivation} There has  been a fruitful interplay between C*-algebras and dynamical systems, which can be
summarized by the diagram
$$
\text{dynamical system $\to$ C*-algebra $\to$ $K$-theoretic invariants.}\label{paradigm}
$$
A classic example is the pair of C*-algebras $F_A$, $O_A$ that
Krieger and Cuntz \cite{CunKri} associated with a shift of finite
type.  They computed the K-groups of these algebras, and observed
that $K_0(F_A)$ was the same as a dimension group that Krieger had
earlier defined in dynamical terms \cite{KriDim}.  Krieger
\cite{KriShift} had shown that this dimension group, with an
associated automorphism, determined the matrix $A$ up to shift
equivalence.  Later it was shown that this dimension group and
automorphism determines the two-sided shift up to eventual
conjugacy, cf. \cite{LinMar}.
 Cuntz \cite{CunMarkovII} showed that the
K-groups of $O_A$ are exactly the Bowen-Franks invariants for flow equivalence. Not only have these algebras
provided interesting invariants for the dynamical systems, but the favor has been returned, in the sense that these
algebras have been an inspiration for many developments in C*-algebras.

Given a piecewise monotonic map $\tau$ of the unit interval, we will associate two C*-algebras $F_\tau$, $O_\tau$, and
calculate their K-theory, thus getting invariants for $\tau$.  The invariants are computable in terms of the
original map and standard concepts from dynamical systems.

\subsection*{Construction}
As our notation indicates, our algebras $F_\tau$ and $O_\tau$ are
generalizations of $F_A$ and $O_A$ respectively.  They include the algebras $F_A$ and $O_A$ as special cases (Corollary \ref{1.47.1}), but also
include simple C*-algebras that are not Cuntz-Krieger algebras (Examples \ref{1.54}, \ref{1.55}).

 Before defining these
algebras, we review the construction of two related algebras $C^*(R(X,\sigma))$ and $C^*(X,\sigma)$.
Suppose that $X$ is a compact Hausdorff space, and $\sigma:X\to X$
is a local homeo\-morphism. Define an equivalence relation
$$R(X,\sigma) = \{(x,y) \in X\times X \mid \sigma^nx = \sigma^n y \hbox{ for some $n \ge 0$}\},$$
 with a suitable
locally compact topology. For $f,g \in C_c(R(X,\sigma))$, define a
convolution product by $(f*g)(x,y) = \sum_{z\sim x} f(x,z)g(z,y)$,
and an involution $*$ by $f^*(x,y) = \overline{f(y,x)}$. Then
$C_c(R(X,\sigma))$ becomes a *-algebra, which can be normed and
completed to become a C*-algebra $C^*(R(X,\sigma))$.  This is a
special case of the  construction of C*-algebras from locally
compact groupoids, so one can make use of the extensive machinery
developed by Renault \cite{RenBook}. The definition of
$C^*(X,\sigma)$ is similar, except that $R(X,\sigma)$ is replaced
by
$$G(X,\sigma) = \{(x,n,y) \in X\times \Z
\times X\mid
\sigma^j x = \sigma^k y \hbox{ and $n = k-j$}\}.$$

This construction is due to Deaconu \cite{Dea}. It generalizes an example of
Renault in \cite{RenBook}, and is generalized further in \cite{Ana, RenCuntzAlg}.   It is  closely
related to a construction of Arzumanian and Vershik \cite{ArzVer} in a measure-theoretic context.  If this construction  is
applied to the local homeo\-morphism given by a one-sided shift of finite type associated with a zero-one matrix $A$, then
$C^*(R(X,\sigma))
\cong F_A$ and
$C^*(X,\sigma)
\cong O_A$, cf.
\cite{Dea,RenCuntzAlg}.

For the construction above to work, it is crucial that $\sigma$ is a local homeo\-morphism. Therefore, if $\tau:I\to I$ is \pwm, we
associate with $\tau$ a local homeo\-morphism $\sigma:X\to X$, constructed by disconnecting the interval $I =[0,1]$ at points in the
forward and backward orbit of the endpoints of intervals of monotonicity. This is a well known technique for one dimensional
dynamical systems, cf. \cite{KeaDisconnect, HofDecomp, RueBook,WalBeta}. Then we define $F_\tau = C^*(R(X,\sigma))$, and
$O_\tau = C^*(X,\sigma)$. (There are alternative, equivalent, ways to define $F_\tau$ and $O_\tau$.  For example, $F_\tau$ can be defined
directly as an inductive limit of interval algebras (cf. \cite[Prop. 12.3]{paper I} and \cite[Thm. 8.1]{DeaSh}), or $F_\tau$ and $O_\tau$
 can be defined as the C*-algebras $\mathcal{F}_E$ and $\mathcal{O}_E$ defined by Pimsner \cite{Pim}, where $E$ is an appropriate
Hilbert module.  (For details, see \cite{DeaMuh}).

 If $\tau$ is surjective, from its construction, $F_\tau$ comes equipped with
a canonical endomorphism
$\Phi$, so we have a C*-dynamical system $(F_\tau,\Phi)$. The algebra $O_\tau$ is isomorphic to the crossed product $F_\tau
\times_\Phi \N$, i.e., the universal C*-algebra containing a copy of $F_\tau$ and an isometry implementing $\Phi$.

\subsection*{Summary of results}

We describe  the
$K$-groups of these algebras in  dynamical terms,  show that they determine the C*-algebras, and show that the
properties of the C*-algebras are closely related to those of the dynamical system $\tau:I\to I$.

\Theo[Props. \ref{1.15}, \ref{1.16}]{$F_\tau$ is always an AI-algebra (Definition \ref{1.16.0}).  It is an AF-algebra iff $\tau$ has no
homtervals.}

(For 
$\tau:I\to I$, an interval
$J$ is a
\emph{homterval} if
$\tau^n$ restricted to
$J$ is a homeo\-morphism for all
$n\ge 0$. The dynamical behavior of an interval map on a homterval is particularly simple, and there are no homtervals if, for example,
$\tau$ is transitive.  Background on homtervals can be found in Collet-Eckmann
\cite[\S II.5]{ColEck} or de Melo-van Strien \cite[Lemma II.3.1]{Melo-S}.)

\Theo[Props. \ref{1.9}, \ref{1.10}]{If $\tau$ is surjective, then $F_\tau$ is simple iff $\tau$ is topologically exact (Definitions \ref{1.7.0},
\ref{1.7.1}), and
$O_\tau$ is simple iff
$\tau$ is transitive.}

If $(X_A,\sigma_A)$ is a transitive one-sided shift of finite type, it is well known that
$X_A$ decomposes into $p<\infty$
pieces permuted cyclically by $\sigma_A$,  such that $\sigma_A^p$ is mixing on each piece. The same kind of
decomposition is valid for transitive
\pwm\ maps (cf. \cite[Cor. 4.7]{paper 2}), and for the associated algebra
$F_\tau$.  However, interestingly, we have to exclude the case where $\tau$ is bijective (as can be seen
from the map $x
\mapsto x+\alpha \bmod 1$, with $\alpha$ irrational).  We say $\tau$ is \emph{essentially injective} if
it is injective on the complement of a finite set of points, which is equivalent to the associated local homeo\-morphism
$\sigma$ being a homeo\-morphism.  Of course, a continuous transitive map of the interval can't be essentially injective.

\Theo[Cor. \ref{1.19}]{If $\tau$ is transitive and not essentially
injective, then $F_\tau$ is a finite direct sum $A_1 \oplus \cdots
\oplus A_n$ of simple AF-algebras, each possessing a unique tracial state. The endomorphism $\Phi$ maps
each $A_i$ to $A_{i-1\bmod n}$.}

\Theo[Thm. \ref{1.35}]{If $\tau$ is transitive and not essentially
injective, then $O_\tau$ is separable, simple, purely infinite,
and nuclear, and it is in the UCT class.}

Under the hypotheses above, by the classification results of Phillips \cite{Phi} and Kirchberg \cite{Kir}, it follows that $O_\tau$
is determined by its K-groups.

Since $F_\tau$ is an AI-algebra, then $K_0(F_\tau)$ is a dimension group, i.e., an inductive limit of ordered groups of the form
$\Z^{n_k}$, and
$K_1(F_\tau) = 0$. If $\tau$ is surjective, the endomorphism $\Phi$ induces an automorphism $\Phi_*$ of $K_0(F_\tau)$. We then refer to
$(K_0(F_\tau), K_0(F_\tau)^+, \Phi_*)$ as the \emph{dimension triple} associated with $\tau$.

In  \cite{paper I, paper 2}, a dimension triple $(DG(\tau), DG(\tau)^+, \L_*)$ is defined in  dynamical
terms, and investigated for various families of interval maps. Here  $\L :C(X,\Z) \to
C(X,\Z)$ is the transfer operator (Definition \ref{1.22}).

\Theo[Cor. \ref{1.29.1}]{If $\tau$ is surjective, then $$(K_0(F_\tau),K_0(F_\tau)^+, \Phi_*) \cong (DG(\tau),DG(\tau)^+, \L_*).$$}

Combining this description of $K_0$, with Paschke's exact sequence \cite{PasKth} for crossed products by endomorphisms, allows
us to describe the K-groups of $O_\tau$.

\Theo[Prop. \ref{1.34}]{If $\tau$ is surjective,
$$
\begin{array}{lclcl}K_0(O_\tau) &\cong  &DG(\tau)/\im(\id-\L_*) &\cong &C(X,\Z)/
\im(\id-\L)\cr
K_1(O_\tau) &\cong&\ker(\id-\L_*) &\cong &\ker(\id-\L).
\end{array}$$}

In the case when $\tau$ is transitive, we can also describe the traces on $F_\tau$ and the KMS-states on $O_\tau$. 

\Theo[Thms. \ref{1.41}, \ref{1.42}]{If $\tau$ is transitive and not essentially injective, there is a unique tracial state on
$F_\tau$ scaled by the canonical endomorphism $\Phi$.  The trace is given by the unique scaled measure on $I$ (Definition \ref{1.36}), and the
scaling factor is $\exp h_\tau$, where $h_\tau$ is the topological entropy of $\tau$. Furthermore, $O_\tau$ has
a unique $\beta$-KMS state for $\beta = h_\tau$, and there is no other $\beta$-KMS state for $0 \le
\beta < \infty$.}

After establishing the results described above, we
 apply them in sections \ref{S13}-\ref{S17} to particular families of maps, including unimodal maps, multimodal maps, interval exchange maps, and
$\beta$-transformations.  
As can be seen from the results above, the K-groups of $F_\tau$ and $O_\tau$ can be recovered from the dynamically
defined dimension triple $(DG(\tau), DG(\tau)^+, \L_*)$. By making use of the study of this dimension triple in  \cite{paper I,
paper 2}, we are able to give explicit descriptions of the K-groups in terms of the map $\tau$ for these families of maps.

\subsection*{Other C*-algebras associated with interval maps}

Other authors have associated C*-algebras
with particular families of interval maps that coincide with ours. For example, for interval exchange maps $\tau:I\to  I$, Putnam
\cite{PutExc1, PutExc2} associated a homeo\-morphism $\sigma:X\to X$, which is identical to the local homeo\-morphism we associated with
$\tau$ viewed as a
\pwm\ map. The crossed product algebra he analyzes is  the same as our $O_\tau$.  

 Katayama-Matsumoto-Watatani
\cite{KatMatWat} defined  C*-algebras $F_\beta^\infty$ and $O_\beta$ for $\beta$-shifts. If $\tau$ is the associated
$\beta$-transformation, for all
$\beta$,
$O_\tau
\cong O_\beta$ (Corollary \ref{1.61}).
 When the orbit of 1 under $\tau$ is finite,
$F_\tau \cong F_\beta^\infty$ (Corollary \ref{1.61}). If the orbit of 1 is infinite, the dimension groups
$K_0(F_\beta^\infty)$ and
$K_0(F_\tau)$ are isomorphic as groups, but it is unknown whether they are order isomorphic, and thus
unknown whether $F_\beta^\infty$ and $F_\tau$ are isomorphic.
The algebras $F_\beta^\infty$ and $O_\beta$ are
 defined in quite different ways than our $F_\tau$ and $O_\tau$; it is only through their K-groups that we are able to establish
these isomorphisms.

Other authors have defined C*-algebras for particular families of interval maps in ways that give different C*-algebras from those
we've defined. Deaconu and
Muhly
\cite{DeaMuh} extended the local homeo\-morphism approach  to branched coverings, which, for example,
includes the full tent map.  The C*-algebra they associate with the
full tent map is not simple; thus it is not the same as our $O_\tau$. Renault
\cite{RenCuntzAlg} investigated groupoids and the associated C*-algebras associated with partially
defined local homeo\-morphisms, which includes
branched coverings as a special case.
Martins-Severino-Ramos \cite{MarSevRam} associated a zero-one matrix $A$ with quadratic maps whose
kneading sequence is periodic, and then computed the K-groups of  $O_A$ in terms of that kneading
sequence. Since the kneading sequence can be periodic without the critical point being periodic, the algebra $O_\tau$ we associate
with such maps is generally different than their
$O_A$, though the computation of $K_0(O_A)$ in terms of the kneading sequence is quite similar.   
For maps on the interval that are open
maps, Exel \cite{Exel} has defined a C*-algebra, which  turns out to be isomorphic to the Cuntz algebra $O_\infty$ for the tent map, while our
$O_\tau$ is isomorphic to $O_2$. Kajiwara and Watatani
\cite{KajWat} have associated a C*-algebra with the complex dynamical systems given by a polynomial, acting on either the whole
Riemann sphere, the Julia set, or its complement.  In a few cases (e.g. $f(z) = z^2-2$) the Julia set is an interval, so that
$f$ on its Julia set can be viewed as an interval map.  For this dynamical system the algebra of Kajiwara and Watatani is isomorphic to
$O_\infty$, while the corresponding algebra $O_\tau$ we define is  isomorphic to $O_2$.

 Another approach that can be taken when $\tau:[0,1]\to [0,1]$ is continuous and surjective
is to form the inverse limit
$(X^\infty,\tau^\infty)$, which will be a homeo\-morphism, and then form the crossed product
$C(X^\infty)\times_{\tau^\infty}\Z$. This is a special case of the approach of Brenken \cite{Bre},
who associates a Cuntz-Pimsner  C*-algebra with closed relations (such as the graph of
a function).  Since the inverse limits will have fixed points, the crossed product algebras
won't be simple, so this also generally gives different algebras than those we consider.

\section{Local homeo\-morphisms associated with piecewise monotonic maps}

In this section we review a construction in \cite{paper I}, which to each piecewise monotonic map $\tau:[0,1]\to
[0,1]$ associates a local homeo\-morphism $\sigma:X\to X$ on $X\subset \R$, where $X$ is constructed by
disconnecting
$[0,1]$ at certain points.

\defi{1.1}{Let
$I = [0,1]$. A map
$\tau:I
\to I$ is  \emph{piecewise monotonic} if there are points
$0 = a_0 < a_1 <
\ldots < a_n = 1$ such that $\tau|_{(a_{i-1},a_i)}$ is continuous and strictly monotonic for $1 \le i
\le n$.  We will assume the sequence $a_0, a_1, \ldots, a_n$ is chosen so that no interval
$(a_{i-1},a_i)$ is contained in a larger open interval on which
$\tau$ is continuous and strictly monotonic. The sequence of points $0 = a_0 < a_1 < \cdots < a_n =
1$ is the {\it partition associated with $\tau$}, and the intervals
$\{(a_{i-1}, a_i)\mid 1\le i\le n\}$ are called the \emph{intervals of monotonicity for $\tau$}. Note
that for
$1
\le i
\le n$, the map
$\tau|_{(a_{i-1},a_i)}$ extends uniquely to a strictly monotonic continuous map $\tau_i:[a_{i-1},a_i] \to I$, which will be
a homeo\-morphism onto its image.}

If $\tau$ is not continuous,  we will ignore the actual values of $\tau$ at the partition
points, and instead view $\tau$ as being multivalued, with the values at $a_i$ (for $1 \le i \le n-1$) being the
values given by left and right limits, i.e. the values of $\tau$ at $a_i$ are $\tau_i(a_i)$ and
$\tau_{i+1}(a_i)$. We define a (possibly multivalued) function $\htau$ on $I$ by setting $\htau(x)$ to
be the set of left and right limits of $\tau$ at $x$.  At points where $\tau$ is continuous, $\htau(x) = \{\tau(x)\}$, and we
identify $\htau(x)$ with $\tau(x)$.
Thus for $A\subset I$, $\htau(A) = \bigcup_i \tau_i(A\cap [a_{i-1},a_i])$, and $\htau^{-1}(A) =
\bigcup_i \{x \in [a_{i-1},a_i] \mid \tau_i(x) \in A\}$.

If $x \in I$,  the \emph{generalized
orbit} of $x$ is the smallest subset of $I$ containing $x$ and forward and backward invariant with respect to
$\htau$. Let
$I_1$ be the union of the generalized orbits of $a_0, a_1,  \ldots, a_n$, and let
$I_0 = I \setminus I_1$.  The set
$I_0$ is dense in
$I$, and both $I_0$ and $I_1$ are forward and backward invariant with respect to $\htau$.

\defi{1.1a}{Let $I = [0,1]$, and  let $I_0$, $I_1$ be as above.
The {\it disconnection of $I$ at points in
$I_1$\/} is the totally ordered set $X$ which consists of  a copy of $I$ with the usual ordering,  but with each point
$x
\in I_1\setminus
\{0,1\}$ replaced by two points $x^- < x^+$. We equip $X$ with the order topology, and define the {\it collapse map\/}
$\pi:X\to I$ by
$\pi(x^\pm) = x$ for $x \in I_1$, and $\pi(x) = x$ for $x \in I_0$. We write  $X_1 = \pi^{-1}(I_1)$, and $X_0 = \pi^{-1}(I_0)
=X\setminus X_1$. }

Here $X$ will be homeomorphic to a compact subset of $\R$.

\prop{1.2}{Let $\tau:I\to I$ be piecewise monotonic, and $X$ as in Definition \ref{1.1a}. Then
\begin{enumerate}
\item $\pi$ is  continuous, and $\pi|_{X_0}$ is a homeo\-morphism from $X_0$ onto $I_0$.\label{I:1.2.1}
\item $X_0$ is dense in
$X$.\label{I:1.2.2}
\end{enumerate}
}

\prooff{\cite[Prop. 2.2]{paper I}}

\notate{If $b_1, b_2 \in I_1$ with $b_1 <
b_2$, then
$I(b_1, b_2)$ is the  order interval
$[b_1^+,b_2^-]_X$. If $b_1 > b_2$, then we define $I(b_1,b_2) = I(b_2,b_1)$, and if $b_1 = b_2$, then we set $I(b_1,b_2) =
\emptyset$. Each set $I(b_1, b_2)$ will be clopen in $X$, and every clopen subset of $X$ is a
finite union of such sets
\cite[Prop. 2.2]{paper I}.}

\prop{1.2.1}{Let $\tau:I\to I$ be a piecewise
monotonic map, with associated partition $a_0 < a_1 < \cdots < a_n$, and let $(X,\pi)$ be as described in
Definition \ref{1.1a}.
\begin{enumerate}
\item  There is a unique continuous map
$\sigma:X\to X$ such that $\pi\circ \sigma = \tau \circ \pi$ on $X_0$.\label{I:1.2.1.1}
\item The sets $X_0$ and $X_1 = X\setminus X_0$ are forward and backward invariant with respect to $\sigma$.\label{I:1.2.1.1.5}
\item $\pi$ is a conjugacy from $\sigma|_{X_0}$ onto $\tau|_{I_0}$.\label{I:1.2.1.2}
\item The sets $J_1 = I(a_0,a_1), \ldots, J_n= I(a_{n-1},a_n)$  are a partition of $X$ into clopen sets such
that for $1\le i\le n$, $\pi(J_i) = [a_{i-1},a_i]$, and  $\sigma|{J_i}$ is a homeo\-morphism from $J_i$
onto
$\sigma(J_i)$.\label{I:1.2.1.3}
\item For $1\le i \le n$, $\pi\circ \sigma = \tau_i\circ \pi$ on $J_i$.\label{I:1.2.1.4}
\end{enumerate}
}

\prooff{\cite[Theorem 2.3]{paper I}}

Note that if $\tau$ is continuous, by (\ref{I:1.2.1.2})  and density of $I_0$ in $I$ and $X_0$ in $X$,  $\pi$ will be a
semi-conjugacy from $(X,\sigma)$ onto $(I,\tau)$.  We will call $\sigma$ the \emph{local homeo\-morphism associated with $\tau$}. The map
$\sigma$ actually has a somewhat stronger property, which will play an important role in our dynamical description of $K_0(F_\tau)$
(Corollary \ref{1.29.1}).

\defi{1.11}{Let $X$ be a compact metric space. A continuous map $\sigma:X\to X$ is a \emph{piecewise
homeo\-morphism} if
$\sigma$ is open and there is a partition of
$X$ into clopen sets, each mapped homeomorphically onto its image.}

If $\tau:I\to I$ is \pwm, and $\sigma:X\to X$ is the associated local homeo\-morphism, then $\sigma$
will be a piecewise homeo\-morphism by Proposition \ref{1.2.1} (\ref{I:1.2.1.3}).

If
$X$ is any compact metric space and is zero dimensional, then any local homeo\-morphism
$\sigma:X\to X$ is a piecewise homeo\-morphism.  (For each point $x$, pick a clopen  neighborhood  on
which $\sigma$ is injective, then select a finite cover of such clopen sets, and finally refine the
cover to form a partition.)

 Property (\ref{I:1.2.1.2})  of Proposition
\ref{1.2.1} can be used to show that $\sigma$ and $\tau$ share many properties. Before being more explicit, we review some
terminology.

\defi{1.7.0}{If $X$ is any topological space, and $f:X\to X$ is
a continuous map, then
$f$ is \emph{transitive} if for each pair
$U, V$ of non-empty open sets, there exists $n \ge 0$ such that
$f^n(U)
\cap V
\not=
\emptyset$. We say $f$ is \emph{strongly transitive} if for
every non-empty open set $U$, there exists $n$ such that
$\cup_{k=0}^nf^k(U) = X$. The map $f$ is \emph{(topologically) mixing } if for every pair $U, V$ of
nonempty open sets, there exists
$N$ such that for all
$n \ge N$,
$f^n(U)\cap V \not= \emptyset$,
 and $f$ is \emph{topologically exact} if for
every non-empty open set $U$, there exists $n$ such that $f^n(U)= X$.}

\defi{1.7.1}{
 If $\tau:I\to I$ is \pwm, we view $\tau$ as undefined at the set $C$
of endpoints of intervals of monotonicity, and say $\tau$ is transitive
if  for each pair
$U, V$ of non-empty open sets, there exists $n \ge 0$ such that
$\tau^n(U)
\cap V
\not=
\emptyset$, and is topologically mixing if  for each pair
$U, V$ of non-empty open sets, there exists $N \ge 0$ such that
for all $n \ge N$, $\tau^n(U)
\cap V
\not=
\emptyset$. We say $\tau$ is {\it strongly transitive\/} if for
every non-empty open set $U$, there exists $n$ such that $\cup_{k=0}^n\htau^k(U) = I$.
 The map $\tau$ is \emph{topologically exact} if
for
every non-empty open set $U$, there exists $n$ such that $\htau^n(U) = I$. (If $\tau$ is continuous and \pwm, these
definitions are consistent with those in Definition \ref{1.7.0}.)}

Clearly every strongly transitive map is transitive, and every topologically exact map is
topologically mixing. The converse implications do not hold, as can be seen by considering the two-sided $n$-shift. However,
every continuous \pwm\ topologically mixing map is topologically exact,  cf. \cite[Thm. 2.5]{Pre}, and
\cite[Thm. 45]{BloCop}.  See also \cite[Prop. 4.9]{paper 2} for an additional result relating mixing and
exactness.

\prop{1.3}{Let $\tau:I\to I$ be \pwm\, and  $\sigma:X\to X$ the associated local homeo\-morphism.
\begin{enumerate}
\item $\sigma$ is surjective iff $\tau$ is surjective.\label{I:1.3.1}
\item $\sigma$ is strongly transitive iff $\tau$ is strongly transitive iff $\tau$ is
transitive iff $\sigma$ is  transitive.\label{I:1.3.2}
\item $\sigma$ is topologically exact iff $\tau$ is topologically exact.\label{I:1.3.3}
\item The topological entropy of $\sigma$ is equal to the topological
entropy of
$\tau$.\label{I:1.3.4}
\end{enumerate}}

\prooff{\cite[Lemma 4.2, Prop. 2.9, Lemma 5.2]{paper I} and \cite[Prop. 2.8]{paper 2}.}

\section{C*-algebras associated with local homeo\-morphisms}

Let $X$ be a compact metric space, and $\sigma:X\to X$ a local homeo\-morphism, i.e., a
continuous open map such that each point admits an open neighborhood on which $\sigma$ is
injective.  In this section we will review the construction of two C*-algebras associated with
$\sigma$. (While this construction is well known, we will need to refer to some of the details later.)

We first describe two locally compact groupoids associated with $\sigma$. (The standard
reference for locally compact groupoids and their associated C*-algebras is the book of Renault
\cite{RenBook}). These groupoids were first described in \cite{RenBook} for the shift map
$\sigma$ on the space
$\Sigma_n$ of sequences
$\{1, 2,
\ldots, n\}^\N$, for
$p$-fold covering maps in \cite{Dea}, for surjective local homeo\-morphisms in \cite{Ana}
and in \cite{ArzRen}, and for general local homeo\-morphisms in \cite{RenCuntzAlg}.
Similar constructions in a measure-theoretic context were introduced in \cite{ArzVer}.

For each
$n
\in
\N$, let
$$R_n = \{(x,y)\in X\times X \mid \sigma^n x = \sigma^n y\},$$
and $R(X,\sigma) = \cup_n R_n.$ Give each $R_n$ the product
topology  from $X \times X$, and give $R(X,\sigma)$ the inductive
limit topology. We write $x\sim y$ when $(x,y) \in R(X,\sigma)$,
and define a product and involution on $C_c(R(X,\sigma))$  by
$$(f\star g)(x,y) = \sum_{z\sim x} f(x,z)g(z,y), \qquad f^*(x,y) = \overline{f(y,x)}.\label{eq1}$$

Each $R_n$ is a compact open subset of $R(X,\sigma)$, so we can
identify $C(R_n)$ with the functions in $C_c(R(X,\sigma))$ whose
support is contained in $R_n$. Thus $C_c(R(X,\sigma)) =
\cup_{n=0}^\infty C(R_n)$. Each subspace $C(R_n)$ will be a
*-subalgebra of $C_c(R(X,\sigma))$. By \cite[Prop.
II.4.2]{RenBook} there is a  norm (necessarily unique) making
$C(R_n)$ a C*-algebra, which we denote by $C^*(R_n)$.  Since
$C_c(R(X,\sigma))$ is the union of the subalgebras $C^*(R_n)$,
there is a unique norm on $C_c(R(X,\sigma))$ satisfying the C*
norm axioms, and the completion is a C*-algebra that we denote
$C^*(R(X,\sigma))$. (By the uniqueness of the C* norm on
$C_c(R(X,\sigma))$, the full and reduced C*-algebras of
$R(X,\sigma)$ coincide; cf. \cite[proof of Prop.
2.4]{RenCuntzAlg}.)  From the discussion above, $C^*(R(X,\sigma))$
is the inductive limit of the algebras $C^*(R_n)$; cf. \cite[Cor.
2.2]{RenRN}.

Now we define a second groupoid
$$G(X,\sigma) = \{(x, n, y) \in X\times \Z \times X \mid \exists\; k, l, \hbox{ with
$\sigma^k x = \sigma^l y$ and $n = k-l$}\},$$
with the product $(x,m,y)(y,n,z) = (x,m+n,z)$,
inverse $(x,n,y)^{-1} = (y,-n, x)$, and with the topology given by the basis of open sets
$$N(U, V, k,l) = \{(x,k-l,y) \in U\times \Z \times V \mid  \sigma^kx =
\sigma^ly\},$$ where  $\sigma^k$ and $\sigma^l$ are injective on the open sets $U$ and $V$. Then we can define a
*-algebra $C_c(G(X,\sigma))$ by a convolution product, and for a suitable norm the completion is a
C*-algebra, denoted by $C^*(G(X,\sigma))$, or  by $C^*(X,\sigma)$.
(The reduced and full C*-algebras associated with  $G(X,\sigma)$
coincide, cf. \cite[Prop. 2.4]{RenCuntzAlg}. We remark that if $\sigma$ is a homeo\-morphism, then the map
$(x,m)\mapsto (x,m,\sigma^mx)$ is a homeo\-morphism and algebraic isomorphism from the  groupoid $X
\times_\phi \Z$ associated with the transformation group
$\Z$ acting on $X$, onto $G(X,\sigma)$. In that context, it follows that $C^*(X,\sigma) \cong C(X) \times_\psi \Z$, where
$\psi(f) = f\circ \sigma^{-1}$, cf. \cite{RenBook}.

There is a close relationship between the two C*-algebras we have
just defined. Suppose that
$\sigma$  is surjective, and let $A = C^*(R(X,\sigma))$.  Then the map $\Phi:C_c(R(X,\sigma)) \to
C_c(R(X,\sigma))$ given by
\begin{equation}\label{(1.1)}
\Phi(f)(x,y) = \frac{1}{\sqrt{p(\sigma(x))p(\sigma(y))}} f(\sigma(x),\sigma(y)),
\end{equation}
where $p(z)$ is the cardinality of $\sigma^{-1}(z)$, extends to a *-endomorphism from $A$ into $A$.
The image of $\Phi$ is
$qAq$, where $q = \Phi(1)$. If $\sigma$ is not injective, then $q \not= 1$. Thus, in the terminology of Paschke
\cite{PasEndo}, $\Phi$ will be a proper corner endomorphism.  Therefore we can  form the crossed product
$A\times_\Phi
\N$ (the universal C*-algebra generated by a copy of $A$ and a non-unitary isometry implementing $\Phi$), and
$C^*(X,\sigma)$ will be isomorphic to $A\times_\Phi
\N$, cf. \cite{Dea} and \cite{Ana}.

Many  maps of the interval that are of interest are not surjective, e.g. the members of
the logistic family
$L_k$ given by $L_k(x)
= kx(1-x)$.  However, these maps are often
\emph{eventually surjective}.  We say a map
$\sigma:X\to X$ is \emph{eventually surjective} if there exists an integer
$n
\ge 0$ such that $\sigma^{n+1}(X) = \sigma^n(X)$. In that case we call $\sigma^n(X)$ the eventual range, and
$\sigma$ restricted to its eventual range is surjective. For example, if $2\le k<4$, the map $L_k$ is not
surjective but is eventually surjective with eventual range being the interval $[0,L_k(1/2)]$.

\label{sect3} If $X$ is a compact Hausdorff space and $\sigma:X\to X$ is an eventually surjective local homeo\-morphism,
with eventual range
$Y$, then the algebra
$C^*(R(X,\sigma))$ is strongly Morita equivalent to $C^*(R(Y,\sigma|_Y))$, and
$C^*(X,\sigma)$ is strongly Morita equivalent to $C^*(Y,\sigma|_Y)$.  This follows from the fact that
each $R(X,\sigma)$ equivalence class meets $Y$, cf.
\cite[Ex. 2.7 and Thm. 2.8]{MuhRenWil}.  Thus for our purposes there is not much lost by restricting
consideration to surjective maps, and we will do that whenever it is  convenient.

\section{Simplicity of $C^*(R(X,\sigma))$ and $C^*(X,\sigma)$.}

 For a groupoid $G$ with  object space $X$, define elements $x, y
\in X$  to be $G$-equivalent if there is an element of $G$ whose source is $x$ and whose range is
$y$.  This is readily seen to be an equivalence relation.  A subset $A$ of $X$ is \emph{$G$-invariant} if $A$ is saturated
with respect to this equivalence relation, i.e., if any element $G$-equivalent to an element of $A$ is itself in $A$.  We
say a locally compact groupoid $G$ is
\emph{minimal} if there are no proper $G$-invariant open subsets.

Let $X$ be a compact Hausdorff space, and $\sigma:X\to X$ a local homeo\-morphism.  Then
$C^*(R(X,\sigma))$ is simple iff $R(X,\sigma)$ is minimal
(\cite[Prop. II.4.6]{RenBook}).  This can be reformulated as follows.

\prop{1.6}{Let $X$ be a compact metric space and $\sigma:X\to X$ a surjective
local homeo\-morphism.  Then
$C^*(R(X,\sigma))$ is simple iff $\sigma$ is topologically exact.}

\prooff{By compactness of $X$,  $R(X,\sigma)$ is minimal iff $\sigma$ is topologically exact, cf. \cite[Prop. 2.1]{KumRen}.}

We remark that surjectivity is not necessary in order for $C^*(R(X,\sigma))$ to be simple.  If $\sigma$ is eventually
surjective, and is topologically exact on its eventual range, then a similar proof shows that $C^*(R(X,\sigma))$ is simple.
Similar modifications are possible for the characterizations of simplicity for $C^*(X,\sigma)$ that follow.

\lem{1.7}{Let $X$ be a compact Hausdorff space, and
$\sigma:X\to X$ a surjective local homeo\-morphism.   The following are equivalent.
\begin{enumerate}
\item For every $x\in X$ and every open set $V$,
there exist
$n,m \in \N$ such that $\sigma^n(x) \in \sigma^m(V)$.\label{I:1.7.1}
\item $G(X,\sigma)$ is minimal.\label{I:1.7.2}
\item $\sigma$ is strongly transitive.\label{I:1.7.3}
\end{enumerate}}

\prooff{(\ref{I:1.7.1}) $\Leftrightarrow$ (\ref{I:1.7.2}) This is clear.

(\ref{I:1.7.1}) $\Rightarrow$ (\ref{I:1.7.3})  By (\ref{I:1.7.1}), for every $x\in X$
and every open set $V$ of $X$, there exist $n, m \in \N$
such that
$x \in \sigma^{-n}\sigma^m(V)$. This is equivalent to
$\bigcup_{m,n} \,\sigma^{-n}\sigma^m(V) = X$.  By
compactness, this holds for a finite subcover, so
$X = \bigcup_{i=1}^k \,\sigma^{-n_i}\sigma^{m_i}(V)$.
Let $M = \max\{n_1, \ldots, n_k\}$.  Then
$$X = \sigma^M(X) = \bigcup_{i=1}^k
\,\sigma^M\sigma^{-n_i}\sigma^{m_i}(V)= \bigcup_{i=1}^k
\sigma^{M-n_i+m_i}(V).$$
Each $M-n_i+m_i \ge 0$, so $\sigma$ is strongly
transitive.

(\ref{I:1.7.3}) $\Rightarrow$ (\ref{I:1.7.1}) Conversely, suppose $\sigma$ is strongly transitive.
Let $V$ be an open  subset of $X$. By strong transitivity, there exists $p$
such that
$X =
\bigcup_{k=0}^p
\sigma^{k}(V)$, so
(\ref{I:1.7.1})  holds.}

Let $X$ be a compact metric space.  A local homeo\-morphism $\sigma:X\to X$  is \emph{essentially free} if for every pair
$m,n$ of distinct non-negative integers there is no open subset of $X$ on which $\sigma^m$ and $ \sigma^n$ agree.
The proof of the next result relies on Renault's result \cite[Prop. 2.5]{RenCuntzAlg} that if
$\sigma:X\to X$ is an essentially free local homeo\-morphism, and $G(X,\sigma)$ is minimal, then
$C^*(X,\sigma)$ is simple.

\prop{1.8}{Let $X$ be a compact metric space, and
$\sigma:X\to X$ a surjective local homeo\-morphism. Then
$C^*(X,\sigma)$ is simple iff  $X$  has
infinite cardinality and $\sigma$ is strongly transitive.}

\prooff{Assume first that $X$ is infinite, and that $\sigma$ is strongly transitive.  Then
$G(X,\sigma)$ is minimal by Lemma \ref{1.7}.  Furthermore, since $\sigma$ is surjective and
transitive,  there exists a point $p$ with a dense forward orbit, cf. \cite[Thm.
5.9]{WalBook}.  If
$\sigma^m =
\sigma^n$ on an open set $V$ with $m < n$, the orbit of $p$ must enter $V$, and thus is
eventually periodic, and so in particular is finite.  The finite orbit of $p$ can't be dense in
the infinite space $X$, so this is impossible.  Thus $\sigma$ is essentially free, and
$G(X,\sigma)$ is minimal, so
$C^*(X,\sigma)$ is simple by  \cite[Prop. 2.5]{RenCuntzAlg}.

Conversely, suppose that $C^*(X,\sigma)$
is simple. Then there are no open $G(X,\sigma)$-invariant subsets by \cite[Prop. II.4.5]{RenBook}, so
$G(X,\sigma)$ is minimal, and therefore  $\sigma$ is strongly transitive by Lemma \ref{1.7}.
Suppose  that $X$  is finite.
Surjectivity of
$\sigma$ implies that $\sigma$ is bijective. Then as observed in the previous section, the map $(x,m)\mapsto
(x,m,\sigma^mx)$ is a homeo\-morphism and isomorphism from the  groupoid $X \times_\sigma \Z$ associated
with the transformation group
$\Z$ acting on $X$, onto $G(X,\sigma)$. Since $\Z$ does not act freely, the associated
C*-algebra $C^*(G(X,\sigma))\cong C(X)\times_\sigma\Z$ is not simple, cf. \cite[p. 230]{Dav}. This
contradiction shows that
$X$ is infinite.}

For an alternative characterization for the case when $\sigma$ is a covering map, see \cite[Thm. 11.2]{ExeVer}.

\section{C*-algebras of piecewise monotonic maps}

\defi{1.9.0}{Let $\tau:I\to I$ be \pwm,  and
let
$\sigma:X\to X$ be the associated local homeo\-morphism. We
define
$F_\tau = C^*(R(X,\sigma))$ and
$O_\tau = C^*(X,\sigma)$.}

 We now give conditions for simplicity of $F_\tau$ and $O_\tau$.

\prop{1.9}{If $\tau:I\to I$ is  \pwm\ and surjective, then $F_\tau$ is simple iff $\tau$ is
topologically exact.}

\prooff{By Proposition \ref{1.3}, $\sigma$ is surjective, and $\tau$ is topologically exact iff $\sigma$ is
topologically exact. Now the result follows from Proposition \ref{1.6}.}

\prop{1.10}{If $\tau:I\to I$ is  \pwm\ and surjective, then $O_\tau$ is
simple iff
$\tau$ is transitive.}

\prooff{By Proposition \ref{1.3},
$\tau$ is surjective, and is  transitive iff
$\sigma$ is strongly transitive. Now the result follows from Proposition \ref{1.8}.}

\section{Structural properties of $F_\tau$}

We are going to show that $F_\tau$ is always an AI-algebra (Definition \ref{1.16.0} below),
describe when
$F_\tau$ is an AF-algebra, and give sufficient conditions for $F_\tau$ to be the direct sum of simple
AF-algebras.

\notate{Let $X$ be a compact metric space, $\sigma:X\to X$, $n \in
\N$, and assume there exists  a partition of $X$ into clopen sets
$X_1, \ldots, X_p\, ,$ and a  clopen set $Y\subset X$, such that
for each $i$, $\sigma^n|_{X_i}$ is a homeo\-morphism from $X_i$ onto
$Y$.  For $1 \le i \le p$, let $\phi_i:Y\to X_i$ be the inverse of
$\sigma^n|_{X_i}$. For each pair of indices $i, j$ we define
$(\phi_i,\phi_j):Y\to R_n(X,\sigma)$ by $(\phi_i,\phi_j)(y) =
(\phi_i(y),\phi_j(y))$. Let $\{e_{ij}\}$ be the standard matrix
units of $M_p = M_p(\C)$. If $Y$ is a compact metric space, we
identify $C(Y,M_n)$ with $M_n(C(Y))$, so that each member $g$ of
$C(Y,M_n)$ can be written as $ g = \sum_{ij} g_{ij} e_{ij}$, where
$g_{ij} \in C(Y)$.}

\prop{1.12}{Assume $X$ is a compact metric space, $\sigma:X\to X$ a local homeo\-morphism, and $n \in
\N$ such that
$X$ admits a partition into clopen sets each mapped homeomorphically onto
$Y
\subset X$ by $\sigma^n$. With the notation above, the map
$\pi:C^*(R_n(X,\sigma))
\to C(Y,M_p)$ given by
$$\pi(f)= \sum_{ij} \left(f\circ (\phi_i,\phi_j)\right)
\  e_{ij}$$ is a *-isomorphism, so that $C^*(R_n) \cong C(Y,M_p)$.}

\prooff{$\pi$ is a *-isomorphism by a straightforward calculation.}

\lem{1.4.0}{If $X$ is a compact metric space, and $\sigma:X\to X$ is a piecewise homeo\-morphism,
then for each $n \in \N$, there exists a partition  $B_1, \ldots, B_q$ of $X$ into clopen sets on
which
$\sigma^n$ is injective, and whose images are pairwise equal or disjoint.}

\prooff{For each point $x \in X$ choose a clopen set on which
$\sigma^n$ is injective, then a finite cover of $X$ by such sets,
and then a partition of $X$ into clopen sets $A_1, \ldots, A_p$ on
which $\sigma^n$ is injective.  Now let $B_1, \ldots, B_q$ be the
partition of $\sigma^nX$ into clopen sets generated by the sets
$\sigma^n(A_i)$ for $1\le i \le p$.  Then  the sets
$\{\sigma^{-n}(B_i)\cap A_j\mid 1\le i \le q,\ 1\le j \le p\}$
form a partition $\P$ of $X$ such that the images of the partition
members under $\sigma^n$ are disjoint or equal.}

\cor{1.13}{Let $X$ be a compact metric space,   $\sigma:X\to X$  a piecewise homeo\-morphism, $n \in
\N$. Choose a partition of $X$ into clopen subsets on which
$\sigma^n$ is injective, with the images of these sets being either equal or disjoint. Let $Y_1,
\ldots, Y_q$ be the collection of the distinct images, and let $n_i$ be the number of inverse images
of each point in $Y_i$. Then $ C^*(R_n)\cong \bigoplus_i C(Y_i,M_{n_i})$.}

\prooff{The sets $\{V_i = \sigma^{-n} Y_i\mid 1\le i \le q \}$ form a partition of $X$ into clopen
$R_n(X,\sigma)$-invariant subsets. Let $R_n|_{V_i}$ denote $R_n \cap (V_i\times V_i)$. Then $C^*(R_n)
\cong
\oplus_i C^*(R_n|_{V_i})$. Now the corollary follows by repeated application of Proposition \ref{1.12},
with $V_i$ in place of $X$ and $Y_i$ in place of $Y$.}

By ``interval" we will always mean a non-degenerate interval, i.e., not a single point.

\lem{1.14}{If $X$ is a compact subset of $\R$, then there is a
decreasing sequence $\{X_n\}$ of subsets of $\R$ such that $\cap_n
X_n = X$, and such that  each $X_n$ is a finite union of closed
intervals. For such a sequence, $C(X)$ is the inductive limit of
the algebras $C(X_n)$, with the connecting maps $C(X_n) \to
C(X_{n+1})$ given by the restriction map.}

\prooff{Let $[a,b]$ be a closed interval containing $X$. Write $\R\setminus X $ as a countable
disjoint union of open intervals $V_1, V_2, \ldots$. Let $X_n =  [a,b] \setminus \bigcup_{i=1}^n
\phi_n(V_i)$, where for each $n$ and each open interval $J$, $\phi_n(J)$ is an open interval with the
same midpoint as $J$, but with length shrunk by a factor $(1-1/n)$. Then
each
$X_n$ is a finite union of closed intervals.  (The shrinking given by $\phi_n$ is needed to insure
that these intervals are non-degenerate). The sets
$X_n$ are a decreasing sequence with intersection equal to $X$. Thus $X$ is the projective limit
$X_1\leftarrow X_2 \cdots$ where the connecting maps are the inclusion  maps. Therefore $C(X)$ is the
inductive limit of the sequence $C(X_1)\rightarrow C(X_2)\rightarrow \cdots$.}

If $\sigma:X\to X$ is the local homeo\-morphism associated with a \pwm\ map $\tau$, then $X$ need not be
totally disconnected.  We recall from \cite{paper I} the characterization of total disconnectedness for
$X$.

\defi{1.15.0}{An interval $J$ is a {\it homterval\/} for a \pwm\ map $\tau:I\to I$ if  $\tau^n$ is a homeo\-morphism
on $J$ for all $n$.}

\exe{1.15.1}{A polynomial of degree $\ge 2$ has homtervals iff it has an attractive periodic orbit, cf.
\cite[Prop. 5.7]{paper I}.}

\prop{1.15}{Let $\tau:I \to I$ be \pwm, and $\sigma:X\to X$ the associated
local homeo\-morphism, and $C$ the associated partition for $\tau$. Then the following are equivalent.
\begin{enumerate}
\item $X$ is totally
disconnected.
\item  The generalized orbit of $C$  is dense in [0,1].
\item $\tau$ has no homtervals.
\end{enumerate}
In particular, these hold if $\tau$ is transitive.}

\prooff{\cite[Prop. 5.8]{paper I}}

 Recall that an AF-algebra is a C*-algebra which is the inductive limit of
a sequence of finite dimensional C*-algebras. The following also can be derived from \cite[Prop.
III.1.15]{RenBook}.

\prop{1.4}{If a compact Hausdorff space $X$ is  totally disconnected, and $\sigma:X\to X$ is a
local homeo\-morphism,  then
$C^*(R(X,\sigma))$ is an AF-algebra.}

\prooff{It suffices to show $C^*(R_n)$ is AF for each $n$. Since
$X$ is totally disconnected, then $\sigma$ is a piecewise
homeo\-morphism.  By Corollary \ref{1.13}, $C^*(R_n)$ is isomorphic
to $\oplus_i C(Y_i,M_n)$, where each $Y_i$ is a clopen subset of
$X$. Since each $Y_i$ is totally disconnected, then each
$C(Y_i,M_n)$ is an AF-algebra, and thus $C^*(R_n)$ and
$C^*(R(X,\sigma))$ are AF-algebras.}

The converse of Proposition \ref{1.4} is false, as was shown by Blackadar \cite[Thm. 7.2.1]{Bla}.

\defi{1.16.0}{A C*-algebra is an \emph{interval algebra} if it is of the form $C([0,1],
A)$  where
$A$ is a finite dimensional C*-algebra, and is an \emph{AI-algebra} if it is an inductive limit of
interval algebras.}

\prop{1.16}{If $\tau:I \to I$ is \pwm, then $F_\tau$ is an AI-algebra. It is an AF-algebra
iff  $\tau$ has no homterval. In
particular, this is true if $\tau$ is  transitive.}

\prooff{First, we show $F_\tau$ is an AI-algebra. Each subalgebra $C^*(R_n)$ is
isomorphic to a sum of algebras of the form $C(Y,M_k)$ for $Y$ a compact subset of $\R$ (Corollary
\ref{1.13}). Here $C(Y)$ is the inductive limit of algebras $C(Y_i)$ where $Y_i$ is a finite
union of closed intervals (Lemma \ref{1.14}), so $C(Y,M_k)$ is the inductive limit of
$C(Y_i,M_k)$, and thus is an AI-algebra.  Hence $\{C^*(R_n)\}$ is an increasing sequence of
AI-algebras, and is dense in $C^*(R(X,\sigma))$, so the latter is an AI-algebra. (An inductive limit of AI-algebras
is an AI-algebra, since AI-algebras admit a local characterization, cf. Elliott's proof of the analogous result for
AT-algebras
\cite[Lemma 4.2 and Thm. 4.3]{Ell2}.)

If $\tau$ has no homtervals, then   $X$ is
totally disconnected (Proposition
\ref{1.15}), and so $F_\tau = C^*(R(X,\sigma))$ is AF
(Proposition \ref{1.4}).

Now suppose that $\tau$ has a homterval $J$.  Let  
 $\pi:X\to I$ be the collapse map. Then $\pi^{-1}(J)\subset X$ is  connected, and $\sigma^n$ is injective on it for all $n$.  The
same is true for the connected component $J'$ of $X$ containing $\pi^{-1}(J)$, and $J'$ is
homeomorphic to $[0,1]$.

Thus the equivalence relation $R(X,\sigma)$ restricted to
$J'$ is trivial.  Let $p$ be the projection corresponding to $J'$  (where $J'$ is viewed as a subset of the
diagonal of $R(X,\sigma)$). Since no distinct points in $J'$ are equivalent, $pC_c(R(X,\sigma))p$
is the set of functions in $C_c(R(X,\sigma))$ with support in the diagonal intersected with
$J'\times J'$. This is then isomorphic to $C(J')$. Since this is dense in
$pC^*(R(X,\sigma))p$, and injections are isometric, then
$pC^*(R(X,\sigma))p\cong C(J')$. If $C^*(R(X,\sigma))$ were an AF-algebra, then the corner
$pC^*(R(X,\sigma))p$ would also be AF. Since
$C(J')$ is not an AF-algebra, then neither is $C^*(R(X,\sigma))$.}

It is also possible to give a constructive proof of the fact that $F_\tau$ is an AI-algebra, by showing that $F_\tau$
is isomorphic to the algebra $A_\tau$ defined in \cite[\S12]{paper I}, cf. \cite[Thm. 8.1]{DeaSh}.

\defi{1.17}{A \pwm\ map $\tau:I\to I$ is \emph{essentially
injective} if there is no pair $U, V$ of disjoint open sets such that $\tau(U) = \tau(V)$.}

The map $\tau$ will be essentially injective iff the associated local homeo\-morphism $\sigma:X\to X$ is injective, cf. \cite[Lemma
11.3]{paper I}.

\theo{1.18}{If $\tau:I\to I$ is transitive and is not essentially injective, and $\sigma:X\to X$ is
the associated local homeo\-morphism, then there is a positive integer
$n$ and  a partition of
$X$ into clopen sets $X_1, \ldots, X_n$ cyclically permuted by $\sigma$, such that $\sigma^n$
restricted to each set
$X_i$ is topologically exact.}

\prooff{\cite[Thm. 4.5]{paper 2}}

\cor{1.19}{Let $\tau:I\to I$ be piecewise monotonic, transitive, and  not essentially
injective, and let $X_1, \ldots, X_n$ be a partition of $X$ as in Theorem \ref{1.18}. Then
$F_\tau= A_1\oplus A_2 \oplus \cdots \oplus A_n$, where for each $i$,
  $A_i \cong C^*(R(X_i,\sigma^n))$ is a simple AF-algebra.
For each index $i$, the canonical endomorphism
$\Phi$ of
$F_\tau$ maps $A_i$ into $A_{i-1 \bmod n}$.}

\prooff{Since $R_k\subset R_{k+1}$ for all $k$, then $R(X,\sigma) = \bigcup_k R_{kn} =
R(X,\sigma^n)$. Since $\sigma$ cyclically permutes $X_1, \ldots, X_n$, then for each $i$,
$\{(x,y) \in R(X,\sigma) \cap(X_i\times X_i)\mid \sigma^n(x)= \sigma^n(y)\}$ is a clopen
invariant subset of
$R(X,\sigma)$, and is isomorphic as a locally compact groupoid to
$R(X_i,\sigma^n)$. For $1 \le i \le n$, let $A_i$ be the closure in $F_\tau$ of the algebra of
functions with compact support in $R(X,\sigma)\cap(X_i\times X_i)$. Then $F_\tau = A_1\oplus
A_2\oplus\cdots\oplus A_n$. Since
$\sigma$ permutes
$X_1,
\ldots, X_n$ cyclically, then
$(x,y)
\in X_i\times X_i$ is in $R(X,\sigma)$ iff
$\sigma^jx = \sigma^jy$ for some $j$ a multiple of $n$. Thus $A_i\cong
C^*(R(X_i,\sigma^n))$ for $1 \le i \le n$.  Since $\sigma^n$ is topologically exact on $X_i$, then each
$A_i$ is a simple AF-algebra  by Proposition \ref{1.6} and Proposition \ref{1.16}.  Since $\sigma$ maps
$X_i$ to $X_{i+1\bmod n}$,  by the
definition of the canonical endomorphism $\Phi$ of $F_\tau$,  when $f\in C_c(R(X,\sigma))$ has
support in $X_i$, then $\Phi(f)$ will have support in
$X_{i-1\bmod n}$. Thus
$\Phi$ will map
$A_i$ into $A_{i-1 \bmod n}$.}

\section{A dynamical description of $K_0(F_\tau)$}

Our goal in this section \label{S7} is to show that $K_0(F_\tau)$ is isomorphic to the dynamically defined
dimension group
$DG(\tau)$ defined in \cite{paper I}.   See
\cite{Goo} or \cite{Eff} for  background on dimension groups.

An ordered abelian group $G$ is \emph{unperforated} if for each positive integer $n$ and $g
\in G$, if $ng \in G^+$ then $g \in G^+$.   An ordered abelian group $G$ satisfies the \emph{Riesz decomposition property} if 
whenever $g_1, g_2, h_1, h_2 \in G$ with $g_i \le h_j$ for
$i, n = 1,2$, then there exists $f$ with $g_i \le f \le h_j$ for $i, j = 1,2$.

\defi{1.21}{An ordered abelian group $G$ is a \emph{dimension group} if
$G$ is unperforated and has the Riesz decomposition property.}

If $G$ is countable, then  by a result of Effros-Handelman-Shen \cite{EffHanShe}, $G$ is a dimension group iff
$G$ is an inductive limit of a sequence of ordered groups $\Z^{n_k}$, cf.  \cite{EffHanShe}.  If $X$ is
any compact Hausdorff space, then
$C(X,\Z)$ (with the pointwise ordering) is easily seen to be a dimension group.

\defi{1.22}{Let $X$ be a compact Hausdorff space and $\sigma:X\to X$ a local homeo\-morphism.  Then the
 \emph{transfer map}
$\L_\sigma:C(X,\Z) \to C(X,\Z)$ is defined by
$$\L_\sigma(f)(x) = \sum_{\sigma y = x} f(y).$$
 (We will write $\L$ in place of $\L_\sigma$ when the meaning is clear from the context.)}

 By compactness of $X$, and  the definition of a local homomorphism, each fiber
$\sigma^{-1}(x)$ is finite, so the sum in the definition of
$\L$ is finite. Note that $\L$ is a positive map, and if $\sigma$ is
injective on a clopen set $E$, then
$$\L_\sigma(\chi_E) = \chi_{\sigma(E)}.$$
 The transfer map will play a central role in our dynamical description of
$K_0(F_\tau)$.

\defi{1.25}{Let $X$ be a compact metric space, and $\sigma:X\to X$ a piecewise homeo\-morphism. Define an
equivalence relation on
$C(X,\Z)$ by
$f\sim g$ if
$\L_\sigma^n f =
\L_\sigma^n g$ for some $n \ge 0$, and denote equivalence class by square brackets. Define $[f]+[g] =
[f+g]$, and $[f] \ge 0$ if $\L_\sigma^n f\ge 0$ for some $n \ge 0$. The set of equivalence classes
is an ordered abelian group, denoted $G_\sigma$.  If $\tau:I\to I$ is \pwm, and $\sigma:X\to X$ is the
associated local homeo\-morphism, then $DG(\tau)$ is defined to be $G_\sigma$. (We will see below that
$G_\sigma$ and $DG(\tau)$ are dimension groups.)}

\lem{1.23}{Let $X$ be a compact Hausdorff space,  $\sigma:X\to X$ a piecewise homeo\-morphism, and $Y = \sigma(X)$. Then
$\L( C(X,\Z)) = \{f\in C(X,\Z) \mid \supp f \subset Y\}$, $\L (C(X,\Z))^+ = \L(C(X,\Z)^+)$, and $\L(C(X,\Z))
\cong C(Y,\Z)$. In particular, $\L (C(X,\Z))$ is a dimension group.}

\prooff{\cite[Lemmas 3.5, 3.9]{paper I}.}

\prop{1.24}{Let $X$ be a compact subset of $\R$, and let $\sigma:X\to X$ be a piecewise
homeo\-morphism. Then
$G_\sigma$ is isomorphic to the inductive limit
$$C(X,\Z)\mapright{\L}\L(C(X,\Z))\mapright{\L}\L^2(C(X,\Z))\mapright{\L}\cdots,$$
and is a dimension group if equipped with the positive cone $G_\sigma^+ = \{[f] \mid f\in
C(X,\Z)^+\}$. The map
$\L_*:G_\sigma\to G_\sigma$ is positive and injective, and is an order automorphism iff
$\sigma$ is eventually surjective.}

\prooff{\cite[Lemma 3.8, Cor. 3.12, and Prop. 4.5]{paper I}}

We are going to show in Proposition \ref{1.27} that $K_0(C(X)) = C(X,\Z)$ for $X$ a compact subset of $\R$.
 We will do this by computing appropriate inductive limits.

\lem{1.26}{If $X_1$, $X_2$, $\ldots$ is a decreasing sequence of
compact subsets of a Hausdorff topological space, and $X =
\bigcap_n\, X_n$, then the
dimension group $C(X,\Z)$ is the inductive limit of the sequence $C(X_1,\Z)\to C(X_2,\Z) \to \cdots \to C(X_n,\Z)\to
\cdots$, where the connecting maps are the restriction maps.}

\prooff{Let $G$ be the inductive limit of the sequence $C(X_n,\Z)$
of dimension groups, thought of as sequences $f_1, f_2, \ldots$
where $f_i\in C(X_i,\Z)$, and such that $f_{i+1} = f_i|_{X_{i+1}}$
for all but finitely many $i$.  Two such sequences are identified
if they agree for all but finitely many terms. If $f_{i+1} =
f_i|_{X_{i+1}}$ for all $i \ge n$, we represent the equivalence
class of this sequence by $[f_n,n]$, which will be the same as
$[g,m]$ if  $f_n|_{X_k} = g|_{X_k}$ for some $k > \max\{m,n\}$.

 Suppose $f \in C(X_n)$ and $f|_X = 0$. Fix $\epsilon > 0$. We
claim
$$\hbox{there exists $k
\ge n$ such that  $|f|< \epsilon$ on $X_k$}.\leqno{(1)}$$ Indeed,
let $G_k = \{x\in X_k
\mid |f(x)| \ge \epsilon\}$; then $G_1, G_2, \ldots$   is a nested decreasing
sequence of compact sets with
empty intersection, so some $G_k$ is empty, which proves (1).

Define $\Psi:G\to C(X,\Z)$ by $\Psi([f,n]) = f|_X$. Suppose that
$\Psi([f,n]) = 0$, i.e., $f \in C(X_n,\Z)$ and $f|_X = 0$.  Choosing
$\epsilon < 1$ in (1), since $f$ has integer values, we conclude that there exists $k>n$ such
that $f|_{X_k} = 0$, so $[f,n] = 0$. Thus $\Psi$
is injective.

To see that $\Psi$ is surjective, let $p$ be the
characteristic function of a clopen set $E \subset X$. Choose $f
\in C(X_1)$ such that $0 \le f \le 1$ and such that $f|_X = p$.
Then $f-f^2 = 0$ on $X$, so by (1), there exists $k \ge 1$ such
that $|f-f^2| \le 3/16$ on $X_k$. Then $f(X_k)$ is contained in
$[0,1/4] \cup [3/4,1]$. Choose a continuous function
$\phi:[0,1]
\to [0,1]$ such that $\phi= 0$ on $[0,1/4]$, and $\phi= 1$ on
$[3/4,1]$. Then $h = (\phi\circ f)|_{X_k}$ is a projection in $C(X_k)$,
and $h|_X = \phi\circ p = p$. Thus $\Psi([h,k])= p$, so $\Psi$ is surjective, which
completes the proof that $C(X,\Z)$ is the desired inductive
limit.}

\prop{1.27}{If $X$ is a compact subset of $\R$, then there is a unique isomorphism  $K_0(C(X)) \cong
C(X,\Z)$, which for each clopen set $E \subset X$ takes the class $[\chi_E]\in K_0(C(X))$ to $\chi_E
\in C(X,\Z)$.}

\prooff{Choose $X_1 \supset X_2 \cdots$ so that each $X_i$ is a finite union of closed intervals, and
such that $\cap_n X_n = X$, cf. Lemma \ref{1.14}. For a closed interval
$J$, we have
$K_0(C(J)) =
\Z$, with a generator being the function constantly 1 on $J$, so we can identify $K_0(C(J))$ with
$C(J,\Z)$.  Therefore
$K_0(C(X_n)) \cong C(X_n,\Z)$, and this isomorphism takes the class $[\chi_E]\in K_0(C(X_n))$ to $\chi_E
\in C(X_n,\Z)$.
 Now the result follows from the fact that $C(X)$ is the inductive
limit of the sequence
$C(X_1)\to C(X_2) \to \ldots$,  continuity of $K_0$ with respect to inductive limits,
 and Lemma \ref{1.26}.}

The second author would like to thank Jack Spielberg for helpful conversations regarding  the proof
of Proposition
\ref{1.27}.

For each clopen subset $E$ of $X$,
we write $\widetilde E $ for the corresponding subset
$\{(x,x) \in X\times X \mid x \in E\}$ of the diagonal of $R(X,\sigma)$, and $[\chi_{\widetilde E}]$
denotes the class of the projection
$\chi_{\widetilde E}$ in
$K_0(C^*(R_n))$.

\lem{1.28}{Let $X$ be a compact subset of $\R$, and let $\sigma:X\to X$ be a piecewise homeo\-morphism. For
each $k \in \N$, there is a unique isomorphism of $K_0(C^*(R_k))$ onto $\L^kC(X,\Z)$ taking
$[\chi_{\widetilde E}]$ to
$\L^k\chi_E$ for each clopen set $E \subset X$, and the elements $[\chi_{\widetilde E}]$ generate
$K_0(C^*(R_k))$.}

\prooff{We may assume without loss of generality that $k = 1$. (Otherwise just replace
$\sigma$ by $\sigma^k$, and observe that $\L_\sigma^k = \L_{\sigma^k}$, and $R_k(X,\sigma) =
R_1(X,\sigma^k)$.)  By Lemma
\ref{1.4.0}, we can construct a partition of
$X$ into clopen sets on each of which $\sigma$ is injective, and whose images are either equal or disjoint. Let $Y_1,
\ldots, Y_p$ be the distinct images of these sets, and let $W_i = \sigma^{-1}(Y_i)$
for $1 \le i \le p$.
By construction,
for each $i$,
$W_i$ admits a finite partition into clopen sets mapped homeomorphically by $\sigma$ onto $Y_i$. It
suffices to prove the statement of the lemma with $W_i$ in place of $X$ for $1\le i\le p$. Thus
without loss of generality, we may assume $p = 1$.

   Let
$X_1,
\ldots, X_n$ be a partition of
$X$ into clopen subsets on which
$\sigma$ is 1-1, such that $\sigma(X_i) = Y \subset X$ for $1 \le i \le n$. We have isomorphisms:
\begin{equation}
K_0(C^*(R_1)) \cong K_0(C(Y,M_n)) \cong K_0(C(Y)) \cong C(Y,\Z) \cong \L C(X,\Z).\label{isos}
\end{equation}
Here the first isomorphism is induced by the isomorphism of $C^*(R_1)$ and $ C(Y,M_n)$ given in
Proposition \ref{1.12}.
The second is induced by the inverse of the isomorphism from $C(Y)$ onto $C(Y)e_{11}$ given by $f \mapsto fe_{11}$. The
third isomorphism is given by Proposition \ref{1.27}. Finally, by Lemma \ref{1.23}, extending functions
in
$C(Y,\Z)$ to functions in $X$ that are zero off $Y$ gives an isomorphism
$C(Y,\Z)
\cong
\L C(X,\Z)$.

Now suppose $E$ is a clopen subset of the partition member $X_1$.  Then under the sequence of
isomorphisms in (\ref{isos}),
$$[\chi_{\widetilde E}] \mapsto [\chi_{\widetilde E} \circ (\phi_1,\phi_1) e_{11}] =
[\chi_{\sigma(E)}e_{11}]
\mapsto
([\chi_{\sigma(E)}] \in C(Y,\Z))
\mapsto \L \chi_E \in \L C(X,\Z).$$
Since the  functions $\L \chi_E$ generate the group $\L C(X,\Z)$, it follows that the classes
$[\chi_{\widetilde E}]$ generate $K_0(C^*(R_1))$, which proves the uniqueness statement in the lemma.}

\theo{1.29}{Let $X$ be a compact subset of $\R$, and let $\sigma:X\to X$ be a piecewise homeo\-morphism.
Then there is a unique isomorphism $K_0(C^*(R(X,\sigma))) \to G_\sigma$ taking
$[\widetilde\chi_E]$ to $ [\chi_E]$ (for
each clopen subset $E$ of $X$).   If
$\sigma$ is surjective, and $\Phi$ is  the canonical endomorphism of
$C^*(R(X,\sigma))$, then the automorphism $\Phi_*$ of $K_0(C^*(R(X,\sigma)))$  is carried
to the automorphism $(\L_\sigma)_*^{-1}$ on $G_\sigma$.}

\prooff{If we identify $X$ with the diagonal of $R(X,\sigma)$, i.e., with
$R_0(X,\sigma)$, for each $n \in \N$ we have the  commutative diagram
$$\begin{array}{ccc}K_0(C^*(R_n)) &\longrightarrow& \L^nC(X,\Z)\cr
\uparrow&&\uparrow{\L^n}\cr
K_0(C^*(R_0)) &\longrightarrow& C(X,\Z)\cr
\end{array},$$
where the vertical arrow on the left is induced by the inclusion of $C^*(R_0)$ in $C^*(R_n)$, and the
horizontal arrows are isomorphisms given in Lemma \ref{1.28}. Since $\L^n:C(X,\Z) \to \L^nC(X,\Z)$  is
surjective, then the map $K_0(C^*(R_0)) \to K_0(C^*(R_n))$ is also surjective. It follows that the
following diagram commutes:
$$\begin{array}{ccc}K_0(C^*(R_{n+1})) &\longrightarrow& \L^{n+1}C(X,\Z)\cr
\uparrow&&\uparrow{\L}\cr
K_0(C^*(R_n)) &\longrightarrow& \L^n C(X,\Z)\cr
\end{array}$$
Since $C^*(R(X,\sigma))$ is the inductive limit of the sequence  $\{C^*(R_n(X,\sigma))\}$, then
$K_0(C^*(R(X,\sigma)))$ is the inductive limit of the sequence
$\{K_0(C^*(R_n))\}$.  By virtue of the commutative diagram above, this is the same as the inductive limit of the sequence
$\L:\L^n C(X,\Z) \to \L^{n+1}C(X,\Z)$, which is in turn isomorphic to $G_\sigma$, cf.
Proposition \ref{1.24}.

Finally, we show that $\L_*\Phi_*$ is the identity map on
$K_0(C^*(R(X,\sigma)))$. For this purpose, it suffices to prove
that $\L_*\Phi_*[\chi_{\widetilde Y}] = [\chi_{\widetilde Y}]$ for
every clopen $Y \subset X$. Without loss of generality, $Y$ is
contained in one of the sets $Y_i$ in Lemma 6.3 (applied with $n =
1$), so that there are clopen sets $E_1, \ldots, E_p$ each mapped
bijectively by $\sigma$ onto $Y$.
 Then
$\Phi(\chi_{\widetilde Y}) = (1/p) \sum \chi_{E_{ij}}$, where
$E_{ij} = \{(x,y) \in E_i \times E_j \mid \sigma x = \sigma y\}$. Letting $w = (1/\sqrt{p}) \sum_i
\chi_{E_{1i}}$, then $ww^* = \chi_{E_{11}}$, and
$$w^*w = (1/p) \sum \chi_{E_{ij}} =\Phi(\chi_{\widetilde Y}).$$
Thus
$$\L_*\Phi_* [\chi_{\widetilde Y}] = \L_* [\chi_{E_{11}}] =  [\chi_{\widetilde
Y}],$$ which completes the proof of the lemma.}

\cor{1.29.1}{If $\tau:I\to I$ is \pwm, then $K_0(F_\tau) \cong DG(\tau)$. If $\tau$ is surjective, and
$\Phi$ is the canonical endomorphism of
$F_\tau$, then the action of $\Phi_*$ on $K_0(F_\tau)$ corresponds to $\L_*^{-1}$ on $DG(\tau)$.}

\prooff{This follows from Theorem \ref{1.29}, and the definitions  $F_\tau= C^*(R(X,\sigma))$ and
$DG(\tau) = G_\sigma$.}

A continuous map $\tau:I\to I$ is \emph{unimodal} if for some $c \in I$, $\tau$ increases on $[0,c]$ and decreases on $[c,1]$. (If
instead $\tau$ decreases and then increases, the map $\phi(x) = 1-x$ is a conjugacy from $\tau$ onto a unimodal map.)

\cor{1.29.2}{If $\tau_i:I\to I$ are unimodal and transitive for $i = 1,2$, and $(K_0(F_{\tau_1}),(\Phi_1)^*) \cong
(K_0(F_{\tau_2}),(\Phi_2)_*)$, then  $\tau_1$ and $\tau_2$ are conjugate.}

\prooff{If $(K_0(F_{\tau_1}),(\Phi_1)^*) \cong
(K_0(F_{\tau_2}),(\Phi_2)_*)$, then by Corollary \ref{1.29.1}, the dimension triples $(DG(\tau_1),DG(\tau_1)^+,\L_*)$ and
$(DG(\tau_2),DG(\tau_2)^+,\L_*)$ are isomorphic, so by \cite[Cor. 8.3]{paper 2}, the maps
$\tau_1$ and $\tau_2$ are conjugate.}

When $\tau$ is transitive, the dimension group $DG(\tau)$ has a unique state scaled by the canonical automorphism $\L_* =
\Phi_*^{-1}$; the scaling factor is $\exp(h_\tau)$, where $h_\tau$ is the topological entropy of $\tau$,
cf. \cite[Thm. 5.3]{paper 2}.
Thus one can recover
$h_\tau$ from the dimension group together with the canonical automorphism.

\cor{1.29.3}{If $\tau_i:I\to I$ are \pwm\ and transitive for $i = 1,2$, and $(K_0(F_{\tau_1}),(\Phi_1)^*) \cong
(K_0(F_{\tau_2}),(\Phi_2)_*)$, then  $h_{\tau_1} = h_{\tau_2}$.}

\prooff{If $(K_0(F_{\tau_1}),(\Phi_1)^*) \cong
(K_0(F_{\tau_2}),(\Phi_2)_*)$, then by Corollary \ref{1.29.1}, the dimension triples $(DG(\tau_1),DG(\tau_1)^+,\L_*)$ and
$(DG(\tau_2),DG(\tau_2)^+,\L_*)$ are isomorphic, so by \cite[Cor. 5.4]{paper 2}, the maps $\tau_1$ and $\tau_2$
have the same topological entropy.}

\section{The module structure of $K_0(F_\tau)$}

If $\tau:[0,1] \to [0,1]$ is \pwm\ and  surjective, then the canonical endomorphism $\Phi$ of
$F_\tau$ induces an automorphism $\Phi_*$ of $K_0(F_\tau)$, cf. Corollary \ref{1.29.1}.
 Under the isomorphism of
$K_0(F_\tau)$ with
$ DG(\tau)$, $\Phi_*^{-1}$ is carried to $\L_*$.

The automorphism $\L_*$  gives
$K_0(F_\tau)\cong DG(\tau)$ the structure of a $\Z[t,t^{-1}]$ module. In \cite{paper I}, a finite set of
generators for this module is described.  For later use, we summarize those
results.

As usual, $\sigma:X\to X$ denotes the local homeo\-morphism associated with $\tau$. If $D$ is a finite
subset of $\R$, we will say distinct points $x, y \in D$ are \emph{adjacent} in
$D$ if there is no element of
$D$ between them. Recall that  $I(c,d) = [c^+,d^-]_X$ for $c, d \in I_1$, where $I_1$ is the
generalized orbit in $I$ of the partition points $a_0, \ldots, a_n$ for $\tau$. We will identify
clopen subsets of $X$ with their characteristic functions, so $I(c,d)$ will be viewed as an element of
$DG(\tau)$. For the definition of $\tau_1, \ldots, \tau_n$ in Theorem \ref{1.31}, see Definition \ref{1.1}.

\theo{1.31}{Let $\tau:I\to I$ be \pwm\ and surjective, with  associated partition $a_0, a_1,
\ldots, a_n$.  Let $M$ be any member of the set $\{\tau_1(a_0)$,
$\tau_1(a_1)$, $\ldots$,
$\tau_n(a_{n-1})$, $\tau_n(a_n)\}$.
Let
$$\J_1 =
\{\,I(c,d)
\mid \hbox{$c$, $d$ are adjacent points in
$\{\,a_0, a_1,
\ldots, a_n, M\,\}$}\,\},$$
and let $\J_2$ be the set of intervals corresponding to jumps at partition points, i.e.,
$$\J_2 = \{\,I(\tau_i(a_i), \tau_{i+1}(a_i)) \mid 1 \le i
\le n-1\,\}.$$ Then
$DG(\tau)$ is generated as a module by $\J_1 \cup \J_2$. }

\prooff{\cite[Thm. 6.2]{paper I}}

\cor{1.32}{If $\tau:I\to I$ is a
continuous, surjective \pwm\ map, with associated partition
$a_0, a_1, \ldots, a_n$, then
$DG(\tau)$ is generated as a module by $I(a_0,a_1)$, $I(a_1,a_2)$, $\ldots$,
$I(a_{n-1},a_n)$.}

\prooff{\cite[Cor. 6.3]{paper I}}

\section{K-groups of $O_\tau$}

 Let
$A$ be a unital C*-algebra,  and $\alpha:A\to A$ a proper corner endomorphism, i.e., $q = \alpha(1)$ is a proper
projection of $A$,  and $\alpha$ is a *-isomorphism from $A$
onto $qAq$.  If  in addition $qAq$ is a full corner of $A$, then there is a cyclic exact sequence due to Paschke
\cite[Thm. 4.1]{PasKth}:
\begin{equation}\label{(1.2)}
 \cdots \to K_j(A) \mapright{\id-\alpha_*} K_j(A) \mapright{i_*} K_j(A\times_\alpha\N) \mapright{}
K_{1-j}(A) \mapright{\id-\alpha_*}\
\cdots \quad (j = 0, 1)
\end{equation}

The following appears in \cite{Dea} for the case where $\sigma$ is a covering map and $G(X,\sigma)$ is minimal.

\prop{1.33.0} {Let $X$ be a compact metric space, and $\sigma:X\to X$ a surjective piecewise homeo\-morphism.  Then
there is an exact sequence
\begin{equation}\label{eq1.33.0}
\begin{matrix}
K_0(C^*(R(X,\sigma))) &\mapright{\id-\Phi_*} &K_0(C^*(R(X,\sigma))) &\mapright{i_*} &K_0(C^*(X,\sigma))\cr
\mapup{} & & & & \mapdown{}\cr
K_1(C^*(X,\sigma)) &\mapleft{i_*} &K_1(C^*(R(X,\sigma))) &\mapleft{\id-\Phi_*}  & K_1(C^*(R(X,\sigma))),
\end{matrix}
\end{equation}
where $\Phi$ is the canonical endomorphism of $C^*(R(X,\sigma))$, cf. equation (\ref{(1.1)}).}

\prooff{If $\sigma$ is injective, then $C^*(R(X,\sigma)) = C(X)$, and $\Phi$ is a *-automorphism of
$C(X)$. Then $C^*(X,\sigma) \cong C(X) \times_\Phi \Z$, and (\ref{eq1.33.0}) is the Pimsner-Voiculescu exact sequence
\cite{PimVoi}. If $\sigma$ is not injective, and $q = \Phi(1)$, then $q \not= 1$, and  $\Phi$ is an isomorphism from
$C^*(R(X,\sigma))$ onto the proper corner $qC^*(R(X,\sigma))q$. Then (\ref{eq1.33.0}) follows from Paschke's exact
sequence (\ref{(1.2)}), if we show that this corner is full, i.e., that the ideal it generates is all of
$C^*(R(X,\sigma))$.   For each $Y$ in Corollary \ref{1.13}, let
$E_Y = \frac{1}{ k(Y)} \sum_{ij}e_{ij}^Y$. Then  $\Phi(1)= \sum_Y E_Y$, and so
 the ideal generated by $q = \Phi(1) \in C^*(R_2)$ is all of $C^*(R_2)$; in
particular that ideal contains $1$. Therefore
$qC^*(R(X,\sigma)) q$ is a full corner of $C^*(R(X,\sigma))$, which completes the proof of the proposition.}

\cor{1.33.2}{If $\tau:I\to I$ is \pwm\ and surjective, there is an
exact sequence
\begin{equation}\label{eq1.33.2}
 \begin{matrix}
K_0(F_\tau) &\mapright{\id-\Phi_*} &K_0(F_\tau) &\mapright{i_*} &K_0(O_\tau)\cr
\mapup{} & & & & \mapdown{}\cr
K_1(O_\tau) &\mapleft{i_*} &K_1(F_\tau) &\mapleft{\id-\Phi_*} &K_1(F_\tau),
\end{matrix}
\end{equation}
where $\Phi$ is the canonical endomorphism of $F_\tau$, cf. equation (\ref{(1.1)}).}

\prooff{This follows from Proposition \ref{1.33.0} and the definitions of $F_\tau$ and $O_\tau$.}

\prop{1.33}{Let $G$ be an ordered abelian group and $T:G\to G$  a
positive homomorphism. Let $G^\infty = \lim\limits_\rightarrow
(T:G\to G)$ be the inductive limit, with canonical homomorphisms
$\mu_k$ from the $k$-th copy of $G$ into $G^\infty$. Then there is
a unique order automorphism $T^\infty$ of $G^\infty$ such that
$T^\infty(\mu_k(z)) = \mu_k(Tz)$ for all $k$, and
\begin{equation}\label{(1.3)}
\coker (\id-T^\infty) \cong
\coker (\id-T),
\end{equation}
\begin{equation}\label{(1.4)}
\ker(\id-T^\infty) \cong \ker (\id-T).
\end{equation} }

\prooff{The proof for the case $G = \Z^n$  in \cite[Prop. 3.1]{CunMarkovII} works for general $G$.}

\prop{1.34}{Let $\tau:I\to I$ be a \pwm\ surjective map, with associated local homeo\-morphism
$\sigma:X\to X$. Let $\Phi$ be the canonical endomorphism of $F_\tau$, and $\L:C(X,\Z) \to C(X,\Z)$
the transfer map.  Then
\begin{equation}\label{(1.5)}
K_0(O_\tau) \cong \coker(\id-\Phi_*) \cong \coker(\id-\L_*)\cong \coker(\id -\L),
\end{equation}
\begin{equation}\label{(1.6)}
K_1(O_\tau) \cong \ker(\id-\Phi_*)\cong \ker (\id-\L_*) \cong \ker (\id - \L).
\end{equation}}

\prooff{We apply the exact sequence (\ref{eq1.33.2}).  Since
$F_\tau$ is an AI-algebra, then
$K_1(F_\tau) = 0$, and the first isomorphisms of (\ref{(1.5)}) and (\ref{(1.6)}) follow. Under the isomorphism
of $K_0(F_\tau)$ with $DG(\tau)$, $\Phi_*$ is carried to $\L_*^{-1}$ (Corollary \ref{1.29.1}), so
$(\id-\Phi_*)K_0(F_\tau)$ is carried to
$$(\id - \L_*^{-1})DG(\tau) = (\id -\L_*)\L_*^{-1}DG(\tau) = (\id -\L_*)DG(\tau).$$
From this, the second isomorphism in (\ref{(1.5)}) follows, and the second isomorphism of (\ref{(1.6)}) follows by
a similar argument.

By Proposition \ref{1.24} and Definition \ref{1.25}, $DG(\tau)$ is isomorphic to the
inductive limit
$\lim\limits_\rightarrow (\L:\L^n(C(X,\Z))
\to
\L^{n+1}(C(X,\Z)))$. Since $\tau$ (and therefore $\sigma$ and $\L$) is surjective, this is the same as
the inductive limit
$\lim\limits_\rightarrow (\L:C(X,\Z) \to
C(X,\Z))$. The induced map $\L^\infty$ on $DG(\tau)$ is $\L_*$, cf. Proposition \ref{1.33}. Thus by
(\ref{(1.3)})
$$
\coker(\id-\L_*) \cong \coker(\id-\L).
$$
This completes the proof of  the third isomorphism in (\ref{(1.5)}).  In a similar manner, from (\ref{(1.4)}) we
get the third isomorphism of (\ref{(1.6)}).}

Applying a result of  Anantharaman-Delaroche \cite{Ana}, Renault \cite[Prop. 2.6]{RenCuntzAlg} has shown that for an
essentially free local homeo\-morphism
$\sigma:X\to X$, if for every non-empty open set $U \subset X$, there exists an open set $V \subset U$, and $m, n \in
\N$ such that  $\sigma^n(V) $ is properly contained in $\sigma^m(V)$, then $C^*(X,\sigma)$ is purely infinite. We
apply this to get a sufficient condition for $O_\tau$ to be purely infinite.

\theo{1.35}{If $\tau:I\to I$ is \pwm, transitive, and not essentially injective,  then
$O_\tau$ is separable, simple, purely infinite, and nuclear, and is in the UCT class $\mathcal{N}$.}

\prooff{Separability follows from the fact that $F_\tau$ is an AF-algebra (Proposition \ref{1.16}), and the fact that
$O_\tau$ is the crossed product of $F_\tau$ by an endomorphism. Simplicity follows from Proposition
\ref{1.10}.  Nuclearity of
$C^*(X,\sigma)$ is established in
\cite{Ana}, or see \cite[Prop. 2.4]{RenCuntzAlg}.  Furthermore, $G(X,\sigma)$ is amenable \cite[Prop.
2.4]{RenCuntzAlg}, so the full and reduced C*-algebras coincide. Since $\tau$ is transitive, then $\sigma$ is
strongly transitive (Proposition \ref{1.3}), so as shown in the proof of Proposition
\ref{1.8},
$\sigma$ is essentially free.  We would like to show  that for each open set $W$ there is an open set $V\subset W$
and positive integers $m, n$ such that
$\sigma^m(V)$ is properly contained in $\sigma^n(V)$. By Theorem \ref{1.18}, there is a partition $X_1,
\ldots, X_p$ of $X$ into clopen sets invariant under $\sigma^p$, such that $\sigma^p$ is topologically
exact on each $X_i$. Choose $X_j$ so that $X_j \cap W \not= \emptyset$, and let $V$ be a proper open subset of
$X_j\cap W$. By exactness of
$\sigma^p$ on
$X_j$, there exists
$k
\ge 0$ such that $(\sigma^p)^kV = X_j \supset V$, which proves that $O_\tau = C^*(X,\sigma)$ is purely infinite.
 Since AI-algebras are in the UCT class $\mathcal{N}$, and that class is
closed under crossed products by $\N$, then $O_\tau$ is in $\mathcal{N}$.}

The hypothesis that $\tau$ is not essentially injective cannot be omitted in Theorem \ref{1.35}. For
example, if $\tau:I\to I$ is the piecewise linear map given by $x \mapsto x+\theta \bmod
1$, with $\theta$ irrational, then $\tau$ is transitive and essentially injective. The
associated local homeo\-morphism $\sigma:X\to X$ is a homeo\-morphism, and $O_\tau \cong C(X)
\times_\Phi
\Z$. Here Lebesgue measure induces a measure on $X$, invariant with respect to $\sigma$,  so
$O_\tau$  has a tracial state. Since $O_\tau$ is simple (Proposition 5.2),
this tracial state  must be faithful, so
$O_\tau$ can't be purely infinite.

For later reference, we recall the classification results of Kirchberg and Phillips (\cite{Kir} and \cite{Phi}).  For
an exposition, see the book of R\o rdam
\cite[\S8.4]{RorBook}.

\theo{1.35.1}{(Kirchberg and Phillips) Let $A_1$ and $A_2$ be unital C*-algebras which are separable, simple, purely infinite, and
nuclear, and are in the UCT class $\mathcal{N}$. If there is a unital isomorphism of $K_0(A_1)$ onto $K_0(A_2)$, and
an isomorphism of
$K_1(A_1)$ onto $K_1(A_2)$, then $A_1$ and $A_2$ are isomorphic.}

By  Theorems  \ref{1.35} and \ref{1.35.1}, when $\tau$ is
\pwm, transitive, and not essentially injective, the C*-algebra $O_\tau$ is determined by
$K_0(O_\tau)$ and $K_1(O_\tau)$.

\section{Traces and KMS-states}

We recall the following definition from \cite{paper 2}.

\defi{1.36}{Let $X$ be a compact metric space, $\psi:X\to X$ a  map that takes Borel sets to Borel sets, and $m$ a
probability measure on
$X$. Then \emph{$\psi$ scales $m$ by a factor $s$} if
$m(\psi(E)) = s\,m(E)$ for all Borel sets $E$ on which
$\psi$ is 1--1.}

\lem{1.37}{Let $\tau:I\to I$ be piecewise monotonic, and $\sigma:X\to X$ the associated local
homeo\-morphism. Let
$m$ be a measure on
$X$ scaled by
$\sigma$ by a factor
$s$. Then there is a  trace on $F_\tau = C^*(R(X,\sigma))$ which satisfies
\begin{equation}\label{(1.8)}
\tr(f) = \int_X
f(x,x)\,dm \hbox{ for
$n \in
\N$ and
$f
\in C^*(R_n)$},
\end{equation}
 and this trace is uniquely determined by (\ref{(1.8)}).}

\prooff{By Lemma \ref{1.4.0}, we can
construct a partition of
$X$ into clopen sets on each of which $\sigma^n$ is injective, and whose images are either equal
or disjoint. Let $Y$ be one of these images, and $E_1, \ldots, E_q$ a partition of $\sigma^{-n}Y$ into
clopen sets mapped bijectively by $\sigma^n$ onto $Y$. Let $\{E_{ij}\}$ be the canonical matrix units
in Proposition \ref{1.12}, i.e., $E_{ij}(x,y) = 1$ if $x \in E_i$ and $y \in E_j$ and
$\sigma^n x = \sigma^n  y$. Define $\tr$ on $\cup_n C^*(R_n)$ by (\ref{(1.8)}).
By Proposition \ref{1.12} and Corollary \ref{1.13}, to prove that $\tr(fg) = tr(gf)$ for all $f, g \in
C^*(R_n)$,  it suffices to show that
$\tr(ef) =
\tr(fe)$ for matrix units
$f, e
\in
\{E_{ij}\}$.  Note that $\tr(E_{ij}) = 0$ if $i\not= j$.  Thus
$\tr(E_{ij}E_{i'j'}) \not= 0$ only if $i = j'$ and $j = i'$, i.e.
$E_{i'j'} = E_{ji}$. Thus we need to verify that $\tr(E_{ij}E_{ji})=
\tr(E_{ji}E_{ij})$, i.e., that $\tr(E_{ii}) = \tr(E_{jj})$. We have
$\tr(E_{ii}) = m(E_i)$, so the key requirement is that $m(E_i) = m(E_j)$ for all
$i, j$.  Since $\sigma^n$ maps $E_i$ and $E_j$ bijectively onto $Y$, and $m$ is scaled by $\sigma$, then $m(Y) =
s^n m(E_i) = s^n m(E_j)$. This completes the proof that $\tr$ defines a trace on each
$C^*(R_n)$. It is clear that $\tr$ is positive and $\tr(1) = 1$, so $\tr$ has norm one on each $C^*(R_n)$, and
therefore extends uniquely to a trace on $C^*(R(X,\sigma))$.}

We now describe one way to find scaling measures.

\defi{1.38}{A map $\tau:I\to I$ is \emph{uniformly piecewise linear} if there is
a number
$s > 0$, and a partition
$0 = a_0 < a_1 <
\cdots a_n = 1$, such that $\tau$ is linear with slope $\pm s$ on each interval $(a_i, a_{i+1})$.}

Note that if $\tau:I\to I$ is uniformly piecewise linear with slopes $\pm s$, then Lebesgue
measure $m$ is scaled by $\tau$ by the factor $s$. Let $\sigma:X\to X$ be the associated local
homeo\-morphism, and define a measure
$\mu$ on
$X$ by
$\mu(A) = m(\pi(A\cap X_0))$. Then it is readily verified that $\mu$ is a probability measure on $X$
scaled by $\sigma$ by a factor $s$, cf. \cite[Prop. 3.3]{paper 2}, and so induces a trace on
$F_\tau= C^*(R(X,\sigma))$ satisfying (\ref{(1.8)}) (with $\mu$ in place of $m$).

Uniformly piecewise linear maps occur more frequently than might be thought.
In the following result,  $\tau$ is assumed to be right continuous at $0$, left continuous at
1, and either left or right continuous at any other points of discontinuity. Recall that the
topological entropy of $\tau$ and $\sigma$ are the same (Proposition \ref{1.3}).

\prop{1.39}{([Parry]) If $\tau:I\to I$ is piecewise monotonic and  transitive,
 then $\tau$ is conjugate to a uniformly piecewise linear map with slopes $\pm
s$, where $s > 1$. The scaling factor $s$ equals $e^{h_\tau}$.}

\prooff{The conjugacy result is in \cite{Par} and \cite[Cor. 4.4]{paper 2}. If $\tau$ is continuous,
then so is the conjugate uniformly piecewise linear map, and the fact that  $\ln s =
h_\tau$, follows from
\cite{MisSzl}. For $\tau$ discontinuous,   see  \cite[Cor. 4.4]{paper 2}.}

\theo{1.41}{Let  $\tau:I\to I$ be piecewise monotonic, not essentially injective, and
 transitive, and let $F_\tau = \oplus A_i$ be the decomposition as a sum of simple
AF-algebras given in Corollary \ref{1.19}. Then each $A_i$ has a unique tracial state  $\tau_i$. There
is a tracial state on $F_\tau$ scaled by $\Phi$, and such a tracial state is unique. The
scaling factor is $e^{-h_\tau }$.}

\prooff{Let $\sigma:X\to X$ be the associated local homeo\-morphism. Let $X_1, \ldots,
X_n$ be a partition of $X$ into clopen sets, cyclically permuted by $\sigma$, such that $\sigma^n$
restricted to each $X_i$ is topologically exact, cf. Theorem \ref{1.18}.

There is a unique state on $DG(\tau) \cong K_0(F_\tau)$ scaled by
$\L_*$ by the factor $s = \exp(h_\tau)$ (\cite[Thm. 5.3]{paper
2}). Since there is a bijection of tracial states on a unital
AF-algebra and states on the associated dimension group, cf.,
e.g., \cite[Prop. 1.5.5]{RorBook}, and since $\Phi_* = \L_*^{-1}$,
it follows that there is a unique tracial state on $F_\tau$ scaled
by $\Phi$ by the factor $1/s = \exp(-h_\tau)$. Furthermore,  by
Theorem \ref{1.29}, $K_0(A_i) =K_0(C^*(R(X_i,\sigma^n)))$ is
isomorphic to $G_{\sigma^n |_{X_i}}$, which has a unique state
(\cite[Thm. 5.3]{paper 2}). It follows that $A_i$ has a unique
tracial state.}

The uniqueness statement regarding the scaled trace in Theorem \ref{1.41} can fail if $\tau$ is
essentially injective. For example, there are minimal interval exchange maps with more than one invariant measure,
cf.
\cite{KeaNonErgodic} or \cite{KeyNew}. Each invariant measure induces a trace on $F_\tau$ scaled by $\Phi$ by the factor $s = 1$.

\cor{1.40}{If $\tau:[0,1]\to [0,1]$ is piecewise monotonic and is topologically exact, then
$F_\tau$ is simple and has a unique tracial state.  The canonical endomorphism $\Phi$ of
 $F_\tau$  scales the trace by a factor $e^{-h_\tau}$, where $h_\tau$ is the
topological entropy of $\tau$.}

\prooff{Topological exactness of $\tau$ implies topological exactness of $\sigma$, so there is just
one summand in the decomposition $F_\tau = \oplus A_i$ in Theorem \ref{1.41}.}

We turn to the question of KMS-states on the algebra $O_\tau$. Suppose  that $\tau:I\to I$ is surjective, so that $O_\tau$ is
isomorphic to
$F_\tau
\times_\Phi
\N$. Then
$O_\tau$ is the universal C*-algebra generated by $F_\tau$ and an isometry $v$ such that $vav^* = \Phi(a)$ for $a \in
F_\tau$.  For each scalar $\lambda$ with $|\lambda| = 1$, there is a unique *-automorphism $\alpha_\lambda$ of
$O_\tau$ which fixes
$F_\tau$ and takes
$v$ to $\lambda v$.  This gives an action of $\T$ on $O_\tau$, or alternatively, an action
of $\R$ (by taking $\lambda = e^{it}$). With respect to this action, the set of fixed points is $F_\tau$, and there is
a faithful conditional expectation from $O_\tau$ onto $F_\tau$ given by $E(b) = \int_{\T}
\alpha_\lambda(b)\,d\lambda$. See
\cite{PasEndo, PasKth}.

\theo{1.42}{If $\tau:I\to I$ is  transitive, and is not essentially injective, then
$O_\tau$  has a unique  $\beta$-KMS state for $\beta = h_\tau $,
 and there is no other $\beta$-KMS state for $0 \le \beta < \infty$. }

\prooff{In the remarks at the end of \cite{PasEndo},  it is observed that $\beta$-KMS states for
$0 < \beta < \infty$ are precisely those of the form $\tr\circ E$, where $\tr$ is a trace scaled by $\Phi$ by the
factor $e^{-\beta}$, and $E$ is the conditional expectation from $O_\tau$ onto $F_\tau$ described above. We know
$F_\tau$ has a unique trace scaled by $\Phi$, and that the scaling factor is $\exp(-h_\tau )$, so the desired result
follows for $0 < \beta$.  For the case $\beta = 0$, a KMS-state is a tracial state on $O_\tau$.  Since the latter is
simple and purely infinite (Theorem \ref{1.35}), no such tracial state exists, so no $\beta$-KMS state exists for $\beta  =
0$.}

If $\sigma$ is exact and positively expansive, the results on KMS-states in Theorem \ref{1.42} follow
from
\cite[Thm. 3.5]{KumRen}.  However, the map $\sigma$ will not be positively expansive
unless the forward orbit under $\tau$ of the set $C$ of endpoints of intervals of monotonicity
is finite, i.e.,  $\tau$ is Markov (as defined below).  (For details, see \cite[remark after Thm. 4.5]{paper 2}).

\section{Markov maps}

\begin{figure}[htb]
\centerline{\includegraphics[scale=0.90]{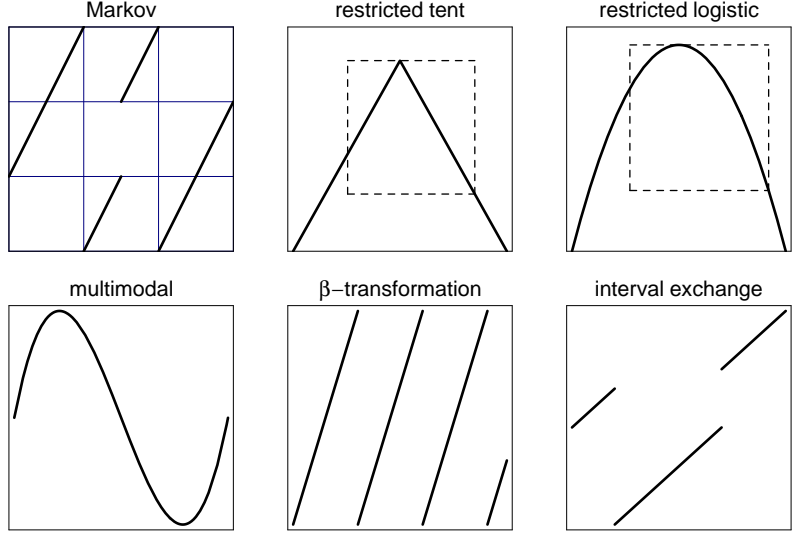}}
\caption{\label{fig1} Examples of piecewise monotonic maps}
\end{figure}

In the \label{S13} remaining sections, we will apply the results developed so far to particular families of interval maps. In
Figure \ref{fig1}, some of the examples that will be discussed are portrayed.

 Recall that if
$\tau:I\to I$ is \pwm, with associated partition $C$, and is discontinuous at $c\in C$,
we view
$\tau$ as multivalued at $c$, with the values given by the left and right limits of $\tau$. For $A \subset [0,1]$, we
write
$\htau(A)$ to denote the set of values of $\tau$ at points in $A$, including the possible multiple values at  points in
$C\cap A$.

\defi{1.43}{Let $\tau:I\to I$ be \pwm. A {\it Markov partition\/} for $\tau$ is
a partition $0 = b_0 < b_1 < \ldots < b_n= 1$, with each $b_i$ being in $I_1$,
such that for each $i$, $\tau$ is monotonic on $(b_i,
b_{i+1})$, such that $\overline{\tau(b_i,b_{i+1})}$ is
 a union of  intervals of the form
$[b_j,b_{j+1}]$, and such that for some $n \ge 0$, $\htau^n(C) \subset \{b_0, b_1, \ldots, b_n\}$. We say $\tau$
is \emph{Markov} if it admits
a Markov partition,  The {\it incidence matrix\/} for $\tau$ with respect to this Markov partition is the zero-one $n\times n$ matrix $A$
where
$A_{ij} = 1$ iff
$\tau(b_{i-1},b_i) \supset (b_{j-1},b_j)$.}

If $\tau$ is \pwm\ with associated partition $C$, and is Markov, then one possible Markov partition is the set of
points in the forward orbit of $C$.  This will be the most common kind of Markov partition that we use, but the slightly
greater generality allowed in the definition above will be useful in showing that all Cuntz-Krieger algebras $O_A$ arise as
$O_\tau$ for some
\pwm\ map $\tau$, cf. Corollary \ref{1.47.1}.

If $\sigma:X\to X$ is the associated local homeo\-morphism, then the
order intervals
$I(b_{i-1},b_i)$ will form a partition of $X$ such that each $\sigma(I(b_{i-1},b_i))$ is a union of some order intervals of
the form
$I(b_{j-1},b_j)$. Since the range of $\tau$ or $\sigma$ will be a union of partition intervals, it follows
that both
$\tau$ and
$\sigma$ will be eventually surjective.

\defi{1.43.1}{If $A$ is an $n \times n$ zero-one  matrix, $G_A$ denotes the stationary inductive limit
 $\Z^n \mapright{A}\Z^n$ in the category of ordered abelian groups. Here $A$ acts by right multiplication,
and
$G_A$ will be a dimension group.  The action of $A$ induces an automorphism of $G_A$ denoted $A_*$.}

\prop{1.44}{If  $\tau:I\to I$ is \pwm\ and Markov, with incidence matrix
$A$,  then $K_0(F_\tau) \cong G_A$. If $\Phi$ is the canonical *-automorphism of $F_\tau$, then the induced
automorphism $\L_* = (\Phi_*)^{-1}$ on $K_0(F_\tau)$ is carried to $A_*$.}

\prooff{This follows from the identification of $K_0(F_\tau)$ with $DG(\tau)$ (Corollary \ref{1.29.1}), and the isomorphism
of $(DG(\tau),\L_*)$ with $(G_A,A_*)$ \cite[Prop. 8.4]{paper I}.}

\cor{1.45}{If  $\tau:I\to I$ is \pwm, surjective, and Markov, with incidence matrix
$A$, then
$$K_0(O_\tau) \cong\coker(\id - A),$$
$$K_1(O_\tau) \cong \ker (\id - A).$$}

\prooff{By Propositions \ref{1.34} and \ref{1.44},
$$K_0(O_\tau) \cong \coker(\id-\L_*) \cong \coker(\id-A_*),$$
which is isomorphic to $\coker(\id-A)$ by Proposition \ref{1.33}. The second isomorphism in the statement of the corollary
follows in a similar fashion.}

If a \pwm\ Markov map $\tau$ is not surjective, then the eventual range $Y$ of $\sigma$ will be a union of some of the
order intervals in the Markov partition.  Restricting
$\sigma$ to its eventual range gives a Markov map whose incidence matrix
$A_Y$ will be the restriction of the original incidence matrix to the rows and columns  corresponding to
intervals in the eventual range. The corresponding algebras $C^*(R(Y,\sigma|_Y))$ and $C^*(Y,\sigma|_Y)$ will
be strongly Morita equivalent to $F_\tau= C^*(R(X,\sigma))$ and $O_\tau= C^*(X,\sigma)$ respectively, as remarked
at the end of Section
\ref{sect3}.  Morita equivalent C*-algebras have the same K-groups, so the isomorphisms in Corollary
\ref{1.45} will hold with $A_Y$ in place of $A$.

Let $X_A$ be the set of
all sequences in $\{1, 2,
\ldots, n\}^\N$ such that $i$ is followed by $j$ only if $A_{ij} = 1$, and let $\sigma_A$ be the
one-sided shift map on $X_A$. We are going to characterize when $(X,\sigma)$ is conjugate to $(X_A,\sigma_A)$.

For each $n \times n$ zero-one matrix $A$ satisfying a certain ``Condition I", Cuntz and Krieger [CK] defined an AF-algebra
$F_A$ and an algebra $O_A$. For finite directed graphs, Condition I is equivalent to the following condition of
Kumjian-Pask-Raeburn \cite[Lemma
3.3]{KumPasRae}. (For a directed graph we will always require that there are no multiple edges between
vertices).

\defi{1.45.1}{A directed graph satisfies \emph{Condition L} if for each cycle $v_1v_2\ldots v_n v_1$, there is an index $i$
and an edge from
$v_i$ to a vertex other than $v_{i+1\bmod n}$.  If $A$ is a \zeroOne\ $n \times n$ matrix, $A$ satisfies Condition L if the
associated directed graph satisfies Condition L.}

\defi{1.46}{Let $\tau:I\to I$ be \pwm\ and Markov, with associated local homeomorphism $\sigma:X\to X$.
Let
$E_1,
\ldots, E_n$ be the associated Markov partition for
$
\sigma:X\to X$, with incidence matrix
$A$. The {\it itinerary map $S:X\to X_A$} is given by $S(x) = s_0s_1s_2\ldots$, where
$\sigma^k(x) \in E_{s_k}$. We say \emph{itineraries separate points of $X$} if the itinerary map
is 1-1. (This is independent of the choice of Markov partition.)}

We now describe when itineraries  separate points. Recall that an interval $J$ is a {\it homterval\/} for a \pwm\ map
$\tau:I\to I$ if  $\tau^n$ is a homeo\-morphism on $J$ for all $n$. A
\pwm\ Markov map
$\tau:I\to I$ is
\emph{piecewise linear} if for each interval $J$ of the Markov partition, there is a partition of $J$ into subintervals on whose
interior
$\tau$ is linear.  We
allow
$\tau$ to be discontinuous at finitely many points in each Markov partition interval, but require that its slope be constant
within each  Markov partition interval.  (See the Markov map in Figure \ref{fig1} on page \pageref{fig1} for an example.)

\prop{1.85}{Let $\tau:I\to I$ be \pwm\ and Markov, with incidence matrix $A$. These are equivalent.
\begin{enumerate}
\item $\tau$ has no homtervals.
\item Itineraries separate points of $X$.
\item $(X,\sigma)$ is conjugate to the one-sided shift of finite type $(X_A,\sigma_A)$.
\end{enumerate}
If $\tau$ is also piecewise linear, these conditions are equivalent to $A$ satisfying Condition L.
}

\prooff{\cite[Props. 8.5 and 8.8]{paper I}}

\prop{1.47}{Let $\tau:I\to I$ be \pwm\ and Markov, and let $\sigma:X\to X$ be the associated local
homeo\-morphism.  If  $A$ is the associated
incidence matrix, and if itineraries  separate points of $X$, then $F_\tau\cong F_A$ and $O_\tau\cong O_A$.}

\prooff{By Proposition \ref{1.85}, $(X,\sigma)$ is conjugate to $(X_A,\sigma_A)$. Thus $C^*(R(X,\sigma)) \cong C^*(R(X_A, \sigma_A))$, and
$C^*(R(X_A,
\sigma_A))\cong F_A$ by (\cite[Example 2]{Dea}). Hence $F_\tau \cong F_A$. In a similar fashion, since $C^*(X_A, \sigma_A)\cong O_A$
\cite[p. 10]{RenCuntzAlg}, then $O_\tau \cong O_A$.}

\cor{1.47.1}{For each $n \times n$ zero-one matrix $A$ satisfying Cuntz-Krieger's Condition I, there is a piecewise linear
Markov map
$\tau$ such that $F_\tau \cong F_A$, and $O_\tau \cong O_A$.}

\prooff{Let $\tau$ be any piecewise linear Markov map, with a Markov partition whose associated incidence matrix is $A$.
Since Condition I is equivalent to Condition L  \cite[Lemma
3.3]{KumPasRae}, the corollary follows from Propositions \ref{1.85} and \ref{1.47}.}

If $\tau$ is piecewise linear, Markov, and surjective, with incidence matrix $A$, we can determine simplicity of $F_\tau$ and
$O_\tau$ from the matrix $A$. Indeed,
 $\tau$ is topologically exact iff $A$ is
primitive (\cite[Cor. 8.9]{paper I}), and this is equivalent to $F_\tau$ being simple (Proposition \ref{1.9}).
Similarly,
$\tau$ is transitive iff
$A$ is irreducible and is not a permutation matrix, and this is equivalent to $O_\tau$ being simple (Proposition
\ref{1.10}).

The next two results will be used later for Markov maps that are unimodal or are
$\beta$-transformations. We write $\Z[t]I(0,1)$ for the set of elements in $C(X,\Z)$ of the form $p(\L)I(0,1)$, with
$p$ a polynomial with integral coefficients.

\lem{1.47.2}{Let $L$ be an endomorphism of an abelian group $G$, and $M$ a subgroup of $G$
invariant under $L$, such that for each $g \in G$ there exists $k \in \N$ such that
$L^k g \in M$. Then
$$\ker ((\id-L)|_M) = \ker (\id-L) \text{ and }\coker ((\id-L)|_M) \cong \coker
(\id-L).$$}

\prooff{Suppose first that $g \in \ker (\id-L)$, so that $Lg = g$.  Choose $k$ so that
$L^k g
\in M$. Then $g = L^k g \in M$, so $g \in \ker ((\id - L)|_M)$. Thus $\ker ((\id-L)|_M) =
\ker (\id-L)$.

Next define $\phi:\coker ((\id-L)|_M) \to \coker(\id-L)$ by
$$\phi(m+ (\id-L)M) = m+ (\id-L)G.$$
If $m \in (\id-L)M$, then $m \in (\id-L)G$, so $\phi$ is well defined. If $m\in M$ and
$\phi(m) = 0$, choose $g \in G$ so that $m = (\id-L)g$. Now choose $k\in \N$ so that
$L^k g \in M$.  Then $L^k m = L^k(g-Lg) = (\id-L)L^k g \in (\id-L)M$. Note that
$$m - L^k m = (\id-L)\sum_{i=0}^{k-1} L^i m \in (\id-L)M.$$
Thus $m = L^k m + (m - L^k m) \in (\id-L)M$, so $\phi$ is injective.
To see that $\phi$ is surjective, let $g \in G$, and again choose $k\in \N$ such that
$L^k g \in M$. Then
$$\phi(L^k g) = L^k g + (\id-L)G = g + (L^k g - g) + (\id-L)G = g + (\id-L)G.$$
Thus $\phi$ is surjective, which completes the proof.}

Recall that a module is \emph{cyclic} if it is generated by a single element, called a cyclic element. If $M$ is a subgroup
of $C(X,\Z)$, invariant under the transfer map $\L$, and $p\in \Z[t]$ is a monic polynomial such that $p(\L)M = 0$, we say
$p$ is the minimal polynomial for $\L$ on $M$ if no polynomial of lower degree annihilates $\L|_M$.  Note that it
may happen that no polynomial annihilates $\L|_M$, so that there is no minimal polynomial for $\L$ on $M$.

\prop{1.47.3}{Let $\tau:I\to I$ be \pwm\ and surjective. Assume that $I(0,1)$ is a cyclic element for the
module
$DG(\tau)$, and that  the minimal polynomial for $\L$ on $M= \Z[t]I(0,1) \subset C(X,\Z)$ is
$$
m(\lambda) = \lambda^p - b_0\lambda^{p-1} - b_1\lambda^{p-2} - \cdots
- b_{p-1}.
$$
Let $n = |m(1)| =  |1-(b_0 + b_1 + \cdots + b_{p-1})|$. \\
If $n \not=0$, then
\begin{equation}
K_0(O_\tau)\cong \Z/n\Z \text{ and } K_1(O_\tau)\cong
0,\label{(1.5.1)}
\end{equation}
and if in addition $\tau$ is transitive, then $O_\tau\cong
O_{n+1}$.\\
If $n = 0$, then
$$K_0(O_\tau) \cong \Z \text{ and } K_1(O_\tau) \cong \Z.$$
Regardless of the value of $n$, $[1]_0$ is a
generator of $K_0(O_\tau)$.}

\prooff{We first show that
\begin{equation}
f \in C(X,\Z) \implies \exists k\in \N \hbox{ such that }\L^k f \in M.\label{(0.9)}
\end{equation}
Since $I(0,1)$ is cyclic for $\L_*$, there exists a Laurent polynomial $q$ with integer
coefficients such that $[f] = q(\L_*)I(0,1)$. Choose $j\in \N$ such that $t^jq(t) \in
\Z[t]$. Then $[\L^j f] = \L_*^j q(\L_*)I(0,1)$, so for some $m \in \N$, $\L^{j+m}f =
\L^{m+j}q(\L)I(0,1)\in M$. Thus (\ref{(0.9)}) holds.

Let $V$ be the  rational linear span of $M$ in
$C(X,\R)$.  Since $v_0 = I(0,1)$ is cyclic for $\L$ on $M$, then for each polynomial
$g \in
\Z[t]$ of degree less than $p$, $g(\L)v_0\not= 0$. (Otherwise $g(\L|_M) = 0$, contrary to the
assumption that
$m$ is the minimal polynomial.)  Therefore, no rational linear combination of
$v_0, \L v_0,
\ldots, \L^{p-1}v_0$ is zero, so these form a basis of $V \cong \Q^p$, and $M$
consists of integral linear combinations of this basis, so $M \cong \Z^p$. For this
basis, the matrix for
$\L|_M$  is
$$B =
\begin{pmatrix}
0  & 1 & 0&\cdots  &0\cr
                0 & 0 & 1&\cdots & 0\cr
                0 & 0 & 0&\cdots & 0\cr
               \cdots    & \cdots &\cdots &\cdots &\cdots  \cr
                0 & 0 & 0&\cdots & 1\cr
                b_0 & b_1& b_2&\cdots & b_{p-1}
\end{pmatrix}
$$
where $B$ acts by right multiplication on $\Q^p$.
Now by elementary operations over $\Z$, (i.e., operations that exchange rows,
multiply a row by -1, or add an integral multiple of one row to another, or  the
analogous column operations), the matrix
$(\id-B)$ can be transformed to a diagonal matrix $D$ with entries $(1, 1, \ldots, 1,
n)$, where $n = |1-(b_0 + b_1 + \cdots + b_{p-1})|$. (This is the Smith form of the
matrix $B$, cf. \cite[p. 248]{LinMar} and \cite{New}).  These
elementary operations can be carried out by left or right multiplication by
invertible integral matrices $W_1$ and $W_2$, with $W_1^{-1}$ and $W_2^{-1}$ also
having integer entries, so that $D = W_1(\id-B)W_2$. It follows that the
kernels of
$D$ and $\id-B$, when viewed as operators on
$\Z^p$ by right multiplication, are isomorphic, and similarly their cokernels are
isomorphic.  Thus if $n \not= 0$, then
$$\coker((\id -\L)|_M) \cong \Z^p/(\Z^p (\id-B))\cong \Z^p/\Z^p D \cong \Z/n\Z$$
$$ \ker((\id-\L)|_M) \cong
\ker(\id-B)\cong \ker D
\cong 0.$$  If $n =0$, then
$$\coker ((\id-\L)|_M)  \cong \Z^p/(\Z^p D) \cong \Z$$
$$\ker((\id-\L)|_M) \cong \ker D\cong \Z.$$ Combining these results with Lemma
\ref{1.47.2} and Proposition \ref{1.34} gives the K-groups described  in the proposition.

Since $M = \Z[t] v_0$,  every element
$v$ of $M$ can be written in the form $p(\L)v_0$, and thus its image in the quotient $M/(\id-\L)M$ is the same as the
image of $p(1)v_0$. Hence  the image of $v_0 = I(0,1)$ generates the quotient.  The element of
$K_0(F_\tau)$ corresponding to the identity in $F_\tau$ is $I(0,1)$, and the inclusion of $F_\tau$ in $O_\tau$ is unital, so
$[1]_0$ generates $K_0(O_\tau)$.

Finally, suppose that $n \not= 0$ and that $\tau$ is  transitive. We would like to apply Theorem \ref{1.35}; for that purpose, we need
to know that
$\tau$ is not essentially injective. If
$\tau$ were essentially injective, then
$\sigma$ would be injective (\cite[Lemma 11.3]{paper I}). Since $\tau$ (and then $\sigma$) are surjective, then
$\sigma$  would be a homeo\-morphism, so $\L I(0,1) = I(0,1)$. Thus the minimal polynomial of $\L$ on $M = \Z[t]
I(0,1)$ would be $m(t) = t-1$,  so $m(1) = 0$, contrary to assumption. Thus
$\tau$ is not essentially injective, so by Theorem \ref{1.35}, $O_\tau$ is simple, separable, nuclear,
purely infinite, and in the UCT class
$\mathcal{N}$, so by the results of Kirchberg \cite{Kir} and Phillips \cite{Phi} (Theorem \ref{1.35.1}), is determined by its
K-groups. By the isomorphisms in (\ref{(1.5.1)}), the K-groups of $O_\tau$ coincide with those of $O_{n+1}$, so these
algebras are isomorphic.}

If $m(t)$ is the minimal polynomial for $\L$ on $M = \Z[t]I(0,1)$,  and $p(t)$ is the minimal polynomial
for the incidence matrix
$A$, then $m$ and $p$ are closely related. In Proposition
\ref{1.47.3},  let
$E_1,
\ldots, E_n$ be the intervals for a Markov partition of
$\sigma$, and
$A$  the incidence matrix. The action of $\L$ on $\sum_i
z_i E_i$ for
$z_i
\in
\Z$ is given by right multiplication by $A$ on $(z_1, z_2,
\ldots, z_n)$.  Since $I(0,1) = E_1 + \cdots + E_n$, then $p(A) = 0$ implies $p(\L)I(0,1) = 0$, so $m(t)$ divides
$p(t)$. On the other hand,
since by hypothesis $I(0,1)$ is a cyclic element for $\L_*$, for each $i$ there exists $f_i \in \Z[t,t^{-1}]$ such
that
$f_i(\L_*)I(0,1) = [E_i]$. For a suitable $k \ge 0$, $f_i(t) t^k \in \Z[t]$, so  for some $K\ge 0$, $\L^K E_i \in
\Z[t]I(0,1)$ for all $i$. Then $m(\L)\L^K E_i = 0$ for all $i$, so $m(A) A^K = 0$. Hence $p(t)$
divides $m(t)t^K$, and thus
 $p(t) = t^j m(t)$ for some
$j\ge 0$. Then $p(1) = m(1)$, so in Proposition \ref{1.47.3}, we could replace $m$ by $p$.

\exe{1.49}{Let $\tau$ be the tent map, i.e., $\tau(x) = 2x$ if $0\le x \le 1/2$, and $\tau(x) = 2-2x$ if $1/2 \le x \le
1$. Then $\tau$ is topologically exact, so $F_\tau$ is a simple AF-algebra and $O_\tau$ is simple and purely infinite.
Furthermore,
$\tau$ is Markov, with incidence matrix
$\begin{pmatrix}1&1\cr1&1\end{pmatrix}$, and the inductive limit of $\Z^2$ with respect to this matrix is the dyadic
rationals,  so  by  Proposition \ref{1.44},
$K_0(F_\tau)$ is isomorphic to the dyadic rationals, with order unit 1, and with the canonical automorphism
being multiplication by 2. Thus $K_0(F_\tau)$ is unitally isomorphic to $K_0(M_{2^\infty})$, and $F_\tau$ and $M_{2^\infty}$ are
unital AF-algebras, so $F_\tau \cong M_{2^\infty}$.
The unique  tracial state of $F_\tau$ is given by
integration with respect to Lebesgue measure, cf. Lemma \ref{1.37}. Since the slopes are $\pm 2$, the entropy of $\tau$ is $\ln
2$ \cite{MisSzl}, so there is a unique
$\beta$-KMS state for $\beta = \ln 2$ (Theorem \ref{1.42}).  By elementary operations over
$\Z$, we can transform
$\id -A$ into the identity matrix, so by Corollary \ref{1.45}, $K_0(O_\tau) \cong \coker(\id) = 0$, and $K_1(O_\tau)
\cong
\ker(\id) = 0$. It follows that $O_\tau  \cong O_A \cong O_2$.  This information on $O_\tau$ can also be derived by viewing
$\tau$ as a unimodal map, cf. Lemma \ref{1.50.1} and Theorem \ref{1.50.5}.}

\exe{1.48}{Let $\tau$ be the Markov map in Figure 1 on page \pageref{fig1}.
Then the associated local homeo\-morphism $\sigma:X\to X$ is Markov with respect to the partition $\{I(0,1/3), I(1/3,
2/3), I(2/3, 1)\}$, with incidence matrix
\begin{equation}\label{(1.10)}
A =
\begin{pmatrix}0&1&1\cr 1&0&1\cr 1&1&0\end{pmatrix}
\end{equation}
All entries of $A^2$ are positive, so $A$ is primitive. By Corollary \ref{1.47.1} and the remarks afterwards,
$F_\tau \cong F_A$ and $O_\tau \cong O_A$, and both $F_\tau$ and $O_\tau$ are simple. By elementary operations over $\Z$, we
can transform $\id -A$ into the diagonal matrix with entries 1, 2, 2, so by Corollary \ref{1.45}, $K_0(O_\tau) \cong
\Z/2\Z
\oplus \Z/2\Z$, and $K_1(O_\tau) \cong 0$.}

Since $K_0(F_\tau) \cong DG(\tau)$, and $DG(\tau)$ is described for most of the following examples in \cite{paper I}
and \cite{paper 2}, in the rest of this paper we will concentrate on describing the algebras $O_\tau$ and their
K-groups.

\section{Unimodal maps}

A map $\tau:[0,1]\to [0,1]$ is \emph{unimodal} if $\tau$ is
continuous, there is a single critical point $c$, and $\tau$
increases on $[0,c]$, and decreases on $[c,1]$. If $\tau$ is
unimodal and surjective, then $\tau(c) = 1$, and we must have
either $\tau(0) = 0$ or $\tau(1) = 0$ (or both). We  will just
consider the case where $\tau(1) = 0$.   (The case where $\tau(0)
= 0$ can be dealt with by similar techniques.)

We first analyze the case where the orbit of $0$ (or equivalently, of the critical point $c$) is periodic.
Let $J_1 = [0,c]$ and $J_{-1} = (c,1]$, and let the itinerary of 0 be the sequence
$n_0,n_1,n_2,\ldots$.  (In other words, $\tau^k 0 \in J_{n_k}$ for $k \ge 0$.) Define $a_k$ to be the product
$n_0n_1\ldots n_k$ for
$k\ge 0$, and $a_{-1} = 1$. With latter convention, note that $a_k a_{k-1} = n_k$ for all $k \ge 0$.

\lem{1.50.1}{Let $\tau:[0,1]\to [0,1]$ be a  surjective unimodal map satisfying
$\tau(1) = 0$. Assume that $0$ is eventually periodic, with $0$, $\tau 0$,
$\ldots$, $\tau^k0$, $\ldots$, $\tau^{p-1}0$ all distinct and
$\tau^{p}0 = \tau^k 0$.  With the notation for the itinerary of 0 described above, the
minimal polynomial $m(t)$ of
$\L$ on  $M = \Z[t]I(0,1)\subset C(X,\Z)$ is as follows.
\begin{enumerate}
\item \label{item50.1} If $0$ is periodic with period $p\ge 3$,
\[
m(t) =t^{p-1} -t^{p-2}-
a_0t^{p-3}- a_1t^{p-4}
-\cdots-a_{p-3}.
\]
\item \label{item50.2} If $0$ is periodic with period 2, $m(t) = t-1$.
\item \label{item50.3} If $0$ is a fixed point of $\tau$, then
$m(t) = t - 2$.
\item \label{item50.4} If $0$ is eventually periodic, with $k > 1$,
\begin{align*}
m(t) = \qquad&t^p -t^{p-1}-
a_0t^{p-2}- a_1t^{p-3}-\cdots-a_{p-2}\cr
 -  \frac{a_{k-1}}{ a_{p-1}}(&t^k -t^{k-1}-a_0t^{k-2}- a_1t^{k-3}
-\cdots-a_{k-2}).
\end{align*}
\item \label{item50.5} If $0$ is eventually periodic, with $k = 1$,
$$m(t) = t^p -t^{p-1}-
a_0t^{p-2}- a_1t^{p-3}-\cdots-a_{p-2}
 -  a_{p-1}(t -1).$$
\end{enumerate}
}

\prooff{Let $c$ be the critical point of $\tau$. Fix $j \ge 0$. If $\tau^j0 \le c$, then
$$\L I( \tau^j0,1 ) = \L I( \tau^j 0,c ) + \L I( c,1 ) =  I( \tau^{j+1}0,1 ) +
I(0,1).$$
If $\tau^j0 > c$, then
$$ \L I( \tau^j0,1 ) =  I( 0,\tau^{j+1}0 ) = I(0,1) -  I( \tau^{j+1}0,1 ).$$
Thus in either case we have
$$\L I( \tau^j0,1 ) = n_j  I( \tau^{j+1}0,1 ) + I(0,1).$$
Rearranging and using $1/n_j = n_j$ gives
\begin{equation} I( \tau^{j+1}0,1 ) = n_j \L I( \tau^j0,1 )- n_jI(0,1).\label{(eqn1)}\end{equation}
By induction,  for $j \ge 1$, $\L$
satisfies
\begin{equation} I( \tau^j0,1 ) = a_{j-1}(\L^j -\L^{j-1}-
a_0\L^{j-2}-a_1\L^{j-3} - \cdots - a_{j-2})I(0,1).\label{(eqn2)}\end{equation}
(When $j = 1$, the right side is  $a_0(\L-1)I(0,1)$.)
Observe that by (\ref{(eqn2)}), for $k \ge 0$, the subgroup generated by
$\{ I( \tau^j0,1 ) \mid 0 \le j\le k\}$ in $C(X,\Z)$ is contained in the subgroup generated by
$\{\L^jI(0,1) \mid 0 \le j\le k\}$.  The opposite containment follows from (\ref{(eqn2)}) by induction. Now we prove
(\ref{item50.1}) -- (\ref{item50.5}).

(\ref{item50.1}) Assume $0$ is periodic with period $p\ge 3$. Then 1 is the unique pre-image of 0, so $\tau^{p-1}0 = 1$.
Putting
$j=p-1$  in (\ref{(eqn2)}) gives
$m(\L)I(0,1) = 0$, where $m$ is as in (\ref{item50.1}).  Here the $p-1$ vectors
$\{ I( \tau^k0,1 )\mid 0 \le k \le p-2\}$ are linearly
independent in $C(X)$, so the vectors $\{\L^kI(0,1)\mid 0 \le k \le p-2\}$ are
independent,
 and thus we
conclude that (\ref{item50.1}) gives the minimal polynomial for
$\L$.

(\ref{item50.2}) If $0$ has period 2, then $\tau 0 = 1$ and $\tau 1 = 0$, so $\L I(0,1) = I(0,1)$. It
follows that $m(t) = t-1$.

(\ref{item50.3}) If $0$ is fixed, then $\L I(0,1) = 2 I( 0,1 )$, so $m(t) =
t-2$ is the minimal polynomial of $\L$.

(\ref{item50.4}) Assume that  0 is eventually periodic, with $k > 1$. Then applying (\ref{(eqn2)}) with $j = p$
and $j = k$, we conclude that
$m(\L)I(0,1) = 0$, where $m$ is as defined in (\ref{item50.4}). Since $I(0,1)$ is a cyclic element for
$\L$ on $M$, then $m(\L) = 0$ on $M$. We show that no polynomial of lower degree
annihilates $\L$ on $M$. The points
$\{\tau^j0\mid 0
\le j
\le p-1\}$ are distinct, and none equals 1,  so the
characteristic functions for the
$p$ intervals
$\{I( \tau^j0,1 )\mid 0 \le j \le p-1\}$  are
linearly independent in $C(X,\R)$. Thus the linear span $V$  of this
set is $p$ dimensional. As remarked above, the linear span of the $p$
vectors
$\{\L^kI(0,1)\mid 0 \le k \le p-1\}$ equals $V$, and thus these vectors also are
linearly independent. Hence, no linear combination of
$1,
\L,
\L^2,
\ldots,
\L^{p-1}$ annihilates $I(0,1)$. It follows that $m$ is the minimal
polynomial of $\L$.

(\ref{item50.5})  Assume that  0 is eventually periodic, with $k = 1$. Then combining (\ref{(eqn2)}) for $j = k=1$
and (\ref{(eqn2)}) for $j = p$ gives the desired formula for $m(t)$.  That $m$ is the minimal polynomial of $\L$ follows in a similar way to
(\ref{item50.4}).}

In the statement of the following theorem, periodic points are also viewed as eventually periodic.

\theo{1.50.5}{Let $\tau:[0,1]\to [0,1]$ be a  surjective unimodal map that satisfies
$\tau(1) = 0$.
\begin{enumerate}
\item\label{1.50.5i} Assume that the critical point $c$ is eventually periodic,
and let $m(t)$ be the minimal polynomial of $\L$ on $\Z[t] I(0,1)$, as described in Lemma \ref{1.50.1}. Let $n = |m(1)|$.
\begin{enumerate}
\item If $n \not=0$, then
$$K_0(O_\tau)\cong \Z/n\Z
\text{ and } K_1(O_\tau)= 0,$$
and if in addition $\tau$ is transitive, then  $O_\tau\cong O_{n+1}$.\\
\item If $n = 0$, then
$$K_0(O_\tau) \cong \Z \text{ and } K_1(O_\tau) \cong \Z.$$
\end{enumerate}
\item \label{1.50.5ii} If 0 is not eventually periodic, then $K_0(O_\tau) \cong \Z$ and $K_1(O_\tau)= 0$. If
in addition $\tau$ is transitive, then $O_\tau \cong O_\infty$.
\end{enumerate}
}

\prooff{By Corollary \ref{1.32}, $DG(\tau)$ is generated as a $\Z[t,t^{-1}]$ module by $I(0,c)$ and $I(c,1)$. Since $\L I(c,1) = I(0,1)$ and $I(0,c)
= I(0,1)-I(c,1)$, then $I(0,1)$ is a cyclic element for $DG(\tau)$.  If
$0$ is eventually periodic, the theorem  follows from  Proposition
\ref{1.47.3} and Lemma \ref{1.50.1}. If $0$ is not eventually periodic, then $K_0(F_\tau) \cong \Z[t,t^{-1}]$ (as abelian groups), with the action of
$\L_*$ given by multiplication by $t$, cf. \cite[Theorem 9.1]{paper I}. Now $K_0(O_\tau) \cong \Z$ and $K_1(O_\tau) \cong 0$
follow from Proposition \ref{1.34}. These are the same K-groups as those for $O_\infty$.
If $\tau$ is transitive, since it is continuous, it can't be essentially injective. Thus  by Theorem \ref{1.35}, $O_\tau$
meets the conditions needed to apply the Kirchberg/Phillips classification results (Theorem \ref{1.35.1}), so $O_\tau \cong O_\infty$.}

Given $1 < s< 2$, we define the {\it restricted tent map\/} $T_s$ by
\begin{equation}
T_s(x) = \begin{cases}
1+s(x-c) & \text{if $x \le c$}\cr
1-s(x-c) &\text{if $x > c$},
\end{cases}
\end{equation}
where $c = 1-1/s$.
(This is the usual symmetric tent map $\tau$ on $[0,1]$ with slopes $\pm s$, restricted to the interval $[\tau^2(1/2), \tau
(1/2)]$, which is the interval of most interest for the dynamics.    Then the map has been rescaled so that its domain
is  [0,1]. See Figure \ref{fig1} on page~\pageref{fig1}.) Note that $T_s(c) = 1$ and $T_s(1) = 0$.

\exe{1.51.0}{Let $\tau= T_s$ be the restricted tent map, with $\sqrt{2} < s < 2$, and assume that the orbit of the critical
point is not eventually periodic. (For example, this occurs if $s = 3/2$, or if $s$ is transcendental.)  By
\cite{JonRan} or
\cite[Lemma 8.1]{paper 2}, $\tau$ is topologically exact. By Theorem \ref{1.50.5}, $O_\tau \cong
O_\infty$.  Thus all such tent maps give isomorphic C*-algebras $O_\tau$.  However, $O_\tau$ comes
equipped with an action of $\mathbb{T}$ (or of $\R$), (see the discussion preceding Theorem \ref{1.42})
and these actions give different KMS-states. In fact, by Theorem
\ref{1.42}, since $\ln s$ is the entropy of $\tau = T_s$, there will be a unique $\beta$-KMS state for $\beta = \ln s$.
Furthermore, the range of the unique state on the dimension group of the AF-algebras
$F_\tau$ will be $\Z[s,s^{-1}]$ (\cite[Prop. 8.2]{paper 2}), as $s$ varies, many of the algebras $F_\tau$
for $\tau = T_s$ will be non-isomorphic, even though the associated algebras $O_\tau$ will all be
isomorphic to $O_\infty$.}

\exe{1.50}{Let $\tau_0(x) = k x(1-x)$, with $c$ the unique critical point, and choose $k \approx 3.68$  so that $p =
\tau_0^3(c)$ is fixed. Let $\tau$ denote $\tau_0|_J$, where  $J = [\tau_0^2(c),
\tau_0(c)]$, rescaled so that the domain of the restricted map is $[0,1]$.  Applying Lemma \ref{1.50.1}~(\ref{item50.5}) (with $k = 1$ and $p = 2$,
since
$\tau^3 c = \tau 0$ is fixed) we conclude that the  minimal polynomial of $\L$ on $\Z[t]I(0,1)$ is $m(t) = t^2-2$. By Theorem \ref{1.50.5},
$K_0(O_\tau) = 0$ and
$K_1(O_\tau) = 0$. By \cite[Example 9.5]{paper I},
$\tau$ is conjugate to $T_{\sqrt{2}}$, which is transitive \cite[Lemma 8.1]{paper 2}, so by Theorem
\ref{1.50.5},
$O_\tau \cong O_2$.}

\section{Multimodal maps}

By a multimodal map we mean a continuous piecewise monotonic map. In this section, we compute the $K$-groups of $O_\tau$ for some multimodal maps.

\exe{1.54}{Let $\tau:I \to I$ be \pwm, continuous, and surjective, with $0 = a_0
< a_1 <
\ldots < a_q = 1$ being the associated partition.  Assume that  the orbits
of
$a_1,
\ldots, a_{q-1}$ are disjoint and infinite, and   $\tau(\{0,1\}) \cap \{0,1\}= \emptyset$.
Then by Corollary \ref{1.29.1} and \cite[Prop. 10.6]{paper I}, $K_0(F_\tau) \cong DG(\tau) \cong  (\Z[t,t^{-1}])^{q-1}$,
with the action of $\L_*$ on $DG(\tau)$ given
by multiplication by
$t$. By Proposition \ref{1.34}, $K_0(O_\tau) \cong \coker(1-t) \cong \Z^{q-1}$, and $K_1(O_\tau) \cong \ker(1-t)
\cong 0$.}

\exe{1.55}{Let $\tau:I\to I$ be continuous, piecewise linear, and topologically mixing,
with slopes
$\pm s$,  with associated partition $a_0 < a_1 < \cdots < a_n$. Assume
\begin{enumerate}
\item $s$ is transcendental,
\item the lengths of the
intervals $[a_0,a_1], \ldots, [a_{n-2},a_{n-1}]$ are independent over $\Z[s,s^{-1}]$,
\item $\tau(\{0,1\}) \cap \{0,1\} = \emptyset$.
\end{enumerate}
\noindent Then by \cite[Prop. 9.2]{paper 2}, $DG(\tau)$ is order isomorphic to $\sum_{i=1}^{n-1}
\Z[s,s^{-1}](a_i-a_{i-1})$,  with the action of $\L_*$ given by multiplication by $s$. As in Example
\ref{1.54}, we find $K_0(O_\tau) \cong \Z^{n-1}$, and $K_1(O_\tau) \cong 0$.}

Observe that for the maps in Examples \ref{1.54} and \ref{1.55},
$K_1(O_\tau)$ is not the torsion free part of $K_0(O_\tau)$, so $O_\tau$ won't be a  Cuntz-Krieger algebra.

\section{Interval exchange maps}

\defi{1.101}{A  \pwm\ map $\tau:[0,1)\to [0,1)$ is an \emph{generalized interval exchange map} if
$\tau$ is increasing on each interval of monotonicity, is  bijective, and is right continuous at
all points.  We will usually identify $\tau$ with its extension to a map from $[0,1]$ into $[0,1]$, defined
to be left continuous at 1.}

Note that if $0 = a_0 < a_1 < \ldots < a_n = 1$ is the partition associated with a generalized interval
exchange map $\tau$, then $\tau$ maps the intervals
$[a_{i-1},a_i)$ to new intervals which partition $[0,1)$.  If a generalized interval exchange map is
linear with slope 1 on each interval of monotonicity, it is called  an \emph{interval exchange map}; in
this case $\tau$ maps each interval $[a_{i-1},a_i)$ to an interval of the same length.

For a generalized interval exchange map, since $\tau:[0,1) \to [0,1)$ is bijective, the associated local
homeo\-morphism $\sigma:X\to X$ will be a homeo\-morphism. Then, as we observed when $C^*(X,\sigma)$ was
defined, $C^*(X,\sigma)$ will be isomorphic to the transformation group C*-algebra $C(X)\times_\psi \Z$,
where the action of $\Z$ is given by the *-automorphism $\psi$ defined by $\psi(f) = f\circ \sigma^{-1}$. If $\tau$ is an interval exchange map, then
our construction of the  homeo\-morphism $\sigma:X\to X$ is the same as that of Putnam \cite{PutExc1}. Thus we
have the following.

\prop{1.56}{If $\tau:I\to I$ is an interval exchange map, then $O_\tau$ coincides with the
transformation group C*-algebra associated with interval exchange maps by Putnam in \cite{PutExc1}.}

When we have referred to the generalized orbit of a point with respect to a \pwm\ map with some
discontinuities, we have previously used the orbit with respect to the multivalued map $\htau$. In the
next proposition, we need to refer to the orbit with respect to $\tau$ itself, and we will call this the
$\tau$-orbit to avoid confusion with our earlier usage.  If a set $B$ has the property that its points
have orbits that are infinite and disjoint, we refer to this property as the IDOC for $B$.

In \cite[Corollary 11.8]{paper I}, the dimension group $DG(\tau)$ was described for generalized interval
exchange maps in which the orbits of the partition endpoints $a_1, a_2, \ldots, a_{n-1}$ satisfy the
IDOC.  Since the associated local homeo\-morphism
$\sigma$ is a  homeo\-morphism, then as a group,
$DG(\tau) = C(X,\Z)$. As a $\Z[t,t^{-1}]$ module, $DG(\tau) \cong \Z[t,t^{-1}]^{n-1}\oplus \Z$, where
the action of $\L_*$ is multiplication by $t$ on the first summand, and trivial on the $\Z$ summand.
We now apply this to compute the K-groups of
$O_\tau$. For interval exchange maps, Proposition \ref{1.57} is due to Putnam \cite{PutExc1}.

\prop{1.57}{Let $\tau:I\to I$ be a generalized interval exchange map, with associated partition
$C= \{a_0, a_1, \ldots, a_n\}$.  If the $\tau$-orbits of $a_1, a_2, \ldots,
a_{n-1}$ are infinite and disjoint, then $K_0(O_\tau) \cong \Z^n$ and $K_1(O_\tau) \cong \Z$.}

\prooff{Let $M = (\Z[t,t^{-1}])^{n-1}$. By the result from \cite[Corollary 11.8]{paper I} described above, $C(X,\Z) =
DG(\tau)\cong M \oplus\Z$, and
$(\id-\L_*)DG(\tau) = (1-t)M \oplus 0$. Thus $C(X,\Z)/(\id-\L)C(X,\Z) \cong M/(1-t)M\oplus \Z$.

Let $\psi:M\to \Z^{n-1}$ be the map $\psi(p_1, \ldots, p_{n-1})=(p_1(1), \ldots,
p_{n-1}(1))$. Then $\psi$ is a homomorphism from $M$ onto $\Z^{n-1}$. The kernel consists of those
$(p_1,
\ldots, p_{n-1}) \in M$ such that $p_i(1) = 0$ for each $i$. If $(p_1, \ldots, p_{n-1}) \in \ker \psi$, for each $i$ we
can write
$p_i(t) = (1-t)q_i(t)$ for $q_i\in \Z[t,t^{-1}]$, so $(p_1, \ldots, p_{n-1}) \in (1-t)M$.
Conversely, each element in $(1-t)M$ is in the kernel of the homomorphism. Thus $M/(1-t)M
\cong \Z^{n-1}$, and so $K_0(O_\tau) \cong \Z^{n-1}\oplus \Z \cong \Z^n$. Finally, $K_1(F_\tau) \cong \ker(1-t)\cong 0\oplus \Z \cong \Z$.}

\section{$\beta$-transformations}

If \label{S17} $\beta > 1$, the {\it $\beta$-transformation\/} on $[0,1)$ is the map $\tau_\beta: x
\mapsto
\beta x
\bmod 1$. (See Figure \ref{fig1} on page~\pageref{fig1}.) We extend $\tau_\beta$ to $[0,1]$ by defining $\tau_\beta(1)= \lim_{x\to 1^-}
\tau_\beta(x)$. Then $\tau_\beta$ is
\pwm, and the associated partition is $\{0, 1/\beta, 2/\beta, \ldots, n/\beta, 1\}$, where $n= [\beta]$ is the greatest integer
$\le \beta$.

Katayama, Matsumoto, and Watatani \cite{KatMatWat} have associated C*-algebras $F_\beta^\infty$ and
$O_\beta$ with the $\beta$-shift. In this section we will show that when $\tau$ is a $\beta$-transformation, the algebras
$O_\tau$ and
$O_\beta$ are isomorphic.  If $1$ is eventually periodic, then $F_\tau \cong F_\beta^\infty$.

The $\beta$-transformation is always topologically exact (\cite[Lemma 10.1]{paper 2}). We need the
following result describing the minimal polynomial of $\L$ on $\Z[t]I(0,1)$.  If $x \in I$, and $n$ is the
greatest integer $\le \beta$, we give the intervals
$[0,1/\beta)$, $[1/\beta,2/\beta)$,
$\ldots$, $[(n-1)/\beta, n/\beta)$, $[n/\beta, 1]$ the labels $0, 1, \ldots, n$, and
define the {\it itinerary\/} of $x\in I$ to be $n_0n_1n_2\ldots$, where
$\tau^kx$ is in the interval with label $n_k$. The calculation of the minimal polynomial in the following result is
similar to our previous calculation for unimodal maps, cf. Lemma \ref{1.50.1}, and so is left to the reader.

\lem{1.58}{(\cite[Lemma 10.2]{paper 2}) Let $\tau$ be the $\beta$-transformation, and $\sigma:X\to X$ the
associated local homeo\-morphism. Assume
$1$ is eventually periodic, with $1, \tau 1, \ldots, \tau^{p-1} 1$
distinct, and with $\tau^p 1 = \tau^k 1$ for
some $k < p$.  Let
the itinerary of 1 be
$n_0 n_1n_2
\ldots n_k \ldots n_{p}\ldots$, and let $M= \Z[t]I(0,1)$. Then $M$ is order isomorphic to $\Z^p$. Let $m(t)$ be the minimal
polynomial for $\L$ restricted to $M$.
\begin{enumerate}
\item If $\tau 1 = 1$, then $m(t) = t-n_0$.
\item If $\tau 1 \not= 1$ and
$\tau^p 1
\not= 0$, then
\begin{equation*}
m(t) = t^p - n_0 t^{p-1} - n_1t^{p-2}
- \cdots - n_{p-1}  - (t^k -
n_0t^{k-1} - \cdots - n_{k-1}).
\end{equation*}
\item If $\tau^p 1 = 0$, then
$$
m(t) = t^{p-1} - n_0 t^{p-2} - n_1t^{p-3}
- \cdots - n_{p-2}.$$
\end{enumerate}
}

As discussed before Theorem \ref{1.42}, there is a natural
action of $\R$ on $O_\tau$, so it makes sense to speak of KMS-states on $O_\tau$. In the following, we use  the
notation of Lemma \ref{1.58}.

\prop{1.60}{Let $\tau$ be the $\beta$-transformation.
\begin{enumerate}
\item If $1$ is eventually periodic,    let
$$n = |m(1)| =
 \begin{cases}
   n_0 - 1 & \hbox{if $\tau 1 = 1$}\cr
   \sum_{i=k}^{p-1} n_i & \hbox{if $\tau 1 \not= 1$ and $\tau^p 1 \not= 0$}\cr
   &\cr
   (\sum_{i=0}^{p-2} n_i)- 1 & \hbox{if $\tau^p 1 = 0$}.
\end{cases}
$$
 Then $n \not= 0$, and  $K_0(O_\tau) \cong \Z/n\Z$ and $K_1(O_\tau) = 0$.
Furthermore, $O_\tau
\cong O_{n+1}$.
\item If $1$ is not eventually periodic, then $K_0(O_\tau) \cong \Z$ and $K_1(O_\tau) = 0$, and $O_\tau \cong
O_\infty$.
\item $F_\tau$ is a simple AF algebra with a unique trace given by Lebesgue measure.
\item The entropy of $\tau$ is $\ln \beta$, and there is a unique $(\ln \beta)$-KMS state for $O_\tau$.
\end{enumerate}}

\prooff{Assume 1 is eventually periodic.  To show $n \not= 0$, we consider three cases. If $\tau 1 = 1$, then $n_0 = \beta \ge 2$, so $m(1) = n_0 - 1
\not= 0$. If
$\tau 1
\not= 1$, and
$\tau^p 1
\not= 0$, suppose that $n= |m(1)| = 0$. Then
$n_k = n_{k+1} =
\cdots = n_{p-1} = 0$. Since
$\sigma^p1 =
\sigma^k 1$, then $n_j = 0$ for all $j \ge k$.  It follows that the orbit of 1 stays in the interval $[0,1/\beta)$
forever. This is possible only if the orbit of 1 lands on 0, which contradicts $\tau^p1 \not= 0$. If $\tau^p1 = 0$,
then $\tau^{p-1}1 = 0$, and $\tau^{p-2}1 =  n_{p-2}/\beta$,  where $n_{p-2} \ge 1$. Since $n_0 \ge 1$, then $n \ge 1$
follows. Thus we've shown that $n \not= 0$. Now (i) follows from Proposition \ref{1.47.3}.

If $1$ is not eventually periodic, then $K_0(F_\tau) \cong DG(\tau) \cong \Z[t,t^{-1}]$, with the action of $\L_*$ given
by multiplication by $t$, cf. \cite[Prop. 10.5]{paper 2}, so the statements about $K_0(O_\tau)$ and
$K_1(O_\tau)$ in (ii) follow from Proposition \ref{1.34}. Since $\tau$ is topologically exact, we can
apply Theorem \ref{1.35} and the Kirchberg-Phillips classification results (Theorem \ref{1.35.1}) to
conclude that $O_\tau \cong O_\infty$.

Since $\tau$ is exact, then (iii) follows from Corollary \ref{1.40}. The fact that the entropy of $\tau$ is $\ln \beta$
can be found in \cite[Prop. 3.7]{paper 2}. The statement about KMS-states follows from Theorem \ref{1.42}.}

If $x \in [0,1]$, then   the
{\it $\beta$-expansion\/} of $x$ is $x = \frac{n_0}{\beta} + \frac{n_1}{\beta^2} + \cdots$,
where
$n_0$ is the greatest integer in
$\beta x$, and for each
$k\ge 1$,
$n_k$ is the greatest integer in $\beta^{k+1} (x-\sum_{i=0}^{k-1} \frac{n_i}{\beta^{i+1}})$.  If $\beta \notin \N$, the
itinerary of $1$  will be the sequence of coefficients $(n_1 n_2\cdots)$ in the $\beta$-expansion of 1.  (If $\beta = n \in \N$, the
itinerary of 1 will  be $(n-1)(n-1)\ldots$, while the $\beta$-expansion of 1 will be $n000\cdots$.) The $\beta$-expansion of 1 will be finite iff
either
$\tau 1 = 1$ or the orbit of 1 lands on 0. Since
$\tau$ is topologically exact, then itineraries separate points, so the
$\beta$-expansion of 1 will be eventually periodic iff 1 is eventually periodic.

\cor{1.61}{If $\tau$ is the $\beta$-transformation, then the C*-algebras  $O_\tau$ and $O_\beta$ are isomorphic. If $1$ is eventually periodic, then
$F_\tau \cong  F_\beta^\infty$.}

\prooff{From Proposition \ref{1.60}, $O_\tau \cong O_{n+1}$ if $1$ is not eventually periodic, and $O_\tau \cong O_\infty$ otherwise. Here $n$ is
described in terms of the itinerary of 1.  In \cite[Cor. 4.13]{KatMatWat}, it also is shown that $O_\beta \cong O_{n+1}$ if the $\beta$-expansion of
$1$ is eventually periodic, and
$O_\beta  \cong O_\infty$ otherwise. A formula for $n$ is given in terms of the $\beta$-expansion of $1$; this formula gives the same value for $n$
as in Proposition \ref{1.60}. Thus $O_\tau \cong O_\beta$ in either case.

If $1$ is eventually periodic, it is shown in \cite[Prop. 10.6]{paper 2} that $DG(\tau)$ and
$K_0(F_\beta^\infty)$ are unitally order isomorphic. Thus $F_\tau$ and $F_\beta^\infty$ are unital
AF-algebras with unitally order isomorphic
$K_0$ groups, so are isomorphic by Elliott's classification of AF-algebras \cite{Ell}.}



\begin{thebibliography}{99}

\bibitem{Ana}
C. Anantharaman-Delaroche, Purely infinite $C\sp *$-algebras arising from dynamical systems. \emph{Bull. Soc. Math.
France} \textbf{125} (1997), no. 2, 199--225.

\bibitem{ArzRen}
V. Arzumanian and J. Renault,  Examples of pseudogroups and their $C\sp *$-algebras. \emph{Operator algebras and quantum
field theory (Rome, 1996)}, 93--104, Internat. Press, Cambridge, MA, 1997.



\bibitem{ArzVer}
V. A. Arzumanian and A. M. Vershik,
Star algebras associated with endomorphisms.
\emph{Operator algebras and group representations}, Vol. I (Neptun,  1980), 17--27,
Monogr. Stud. Math., \textbf{17},
Pitman, Boston, MA, 1984.

\bibitem{Bla} B. Blackadar,
 Symmetries of the CAR algebra.
 \emph{Ann. of Math.} (2)  \textbf{131}  (1990),  no. 3, 589--623.


\bibitem{Bre}
B. Brenken, C*-algebras associated with topological relations, preprint.

\bibitem{BloCop}
L. S. Block and W. A. Coppel,
\emph{Dynamics in one dimension.}
Lecture Notes in Mathematics, \textbf{1513}.
Springer-Verlag, Berlin, 1992.

\bibitem{BoyFieFie}
M. Boyle, D. Fiebig, and U. Fiebig, A dimension group for local
homeo\-morphisms and endomorphisms of one-sided shifts of finite type.  \emph{J. Reine Angew.
Math.}  \textbf{487} (1997), 27--59.

\bibitem{ColEck}
P. Collet and J. Eckmann, \emph{Iterated maps on the interval as dynamical systems}, Progress in Physics 1, Birkh\" auser,
Boston, 1980.

\bibitem{CunMarkovII}
J. Cuntz,  A class of $C\sp{*} $-algebras and topological Markov chains. II.  Reducible chains and the Ext-functor for
$C\sp{*} $-algebras. \emph{Invent. Math.} \textbf{63} (1981),  no. 1, 25--40.

\bibitem{CunKri}
J. Cuntz and W. Krieger, A class of C*-algebras and topological Markov chains.
\emph{Invent. Math.}  \textbf{56}  (1980), no. 3, 251--268.


\bibitem{Dav} K. Davidson,  \emph{C*-algebras by example.} Fields Institute Monographs, 6.
American Mathematical Society, Providence, RI, 1996.


\bibitem{Dea}
V. Deaconu,
Groupoids associated with endomorphisms.
\emph{Trans. Amer. Math. Soc.} \textbf{347} (1995), no. 5, 1779--1786.

\bibitem{DeaMuh} V. Deaconu and P. Muhly, C*-algebras associated with branched coverings.  \emph{Proc. Amer. Math. Soc.}  \textbf{129}
(2001),  no. 4, 1077--1086 (electronic).

\bibitem{DeaSh} V. Deaconu and F. Shultz, C*-algebras associated with interval maps, 40 pages,  \emph{arXiv:math.} OA/0405469.

\bibitem{Ur-De}
M. Denker,  and M. Urba\'nski, On the existence of conformal measures.  \emph{Trans. Amer. Math.
Soc.}  \textbf{328}  (1991),  no. 2, 563--587.

\bibitem{Eff}   E. Effros, \emph{Dimensions and C*-algebras.} CBMS Regional Conference Series in
Mathematics, 46. Conference Board of the Mathematical Sciences, Washington, D.C., 1981.

\bibitem{EffHanShe}
E. G. Effros, D. E. Handelman, and C. L. Shen, Dimension groups
and their affine representations. \emph{Amer. J. Math.}
\textbf{102} (1980),  no. 2, 385--407.

\bibitem{Ell}   G. Elliott, On the classification of inductive limits of sequences of semisimple
finite-dimensional algebras.  \emph{J. Algebra}  \textbf{38}  (1976), no. 1, 29--44.

\bibitem{Ell2}   G. Elliott, On the classification of C*-algebras of real rank zero, \emph{J. Reine Angew. Math.} \textbf{443} (1993), 179-219.

\bibitem{Exel} R. Exel, personal communication.

\bibitem{ExeVer} R. Exel and A. Vershik, C*-algebras of irreversible dynamical systems.
\emph{arXiv:math.} OA/0203185 v1 10 Mar 2002.

\bibitem{GioPutSka} T. Giordano, I. Putnam, and C. Skau, Topological orbit equivalence and C*-crossed products.
 J. Reine Angew. Math. 469  (1995), 51--111.


\bibitem{Goo}
K. R. Goodearl, \emph{Partially ordered abelian groups with interpolation.} Mathematical Surveys and Monographs,
20. American Mathematical Society, Providence, RI, 1986.


\bibitem{HofDecomp}
    F. Hofbauer, The structure of piecewise monotonic transformations.  \emph{Ergodic Theory
Dynamical Systems}  \textbf{1}  (1981), no. 2, 159--178.

\bibitem{JonRan}
L. Jonker and D. Rand, Bifurcations in one dimension. I. The nonwandering set.  \emph{Invent.
Math.}  \textbf{62}  (1981), no. 3, 347--365.

\bibitem{KajWat} T. Kajiwara and Y. Watatani, C*-algebras associated with complex dynamical systems.
\emph{arXiv:math.OA/0309293}.

\bibitem{KatMatWat}
 Y. Katayama,  K. Matsumoto, and  Y. Watatani,
 Simple $C\sp *$-algebras arising from $\beta$-expansion of real numbers.
  \emph{Ergodic Theory Dynam. Systems}  \textbf{18}  (1998),  no. 4, 937--962.

\bibitem{Kat1} T. Katsura,  A class of C*-algebras generalizing both graph algebras and homeo\-morphism C*-algebras I,
fundamental results. 37 pages,  \emph{arXiv:math.} OA/0207252.

\bibitem{Kat2} T. Katsura,  On C*-algebras associated with C*-correspondences,  29 pages,
\emph{arXiv:math.} OA/0309088.

\bibitem{KeaDisconnect}
M. Keane, Interval exchange transformations.  \emph{Math. Z.}  \textbf{141} (1975), 25--31.

\bibitem{KeaNonErgodic}
M. Keane,  Non-ergodic interval exchange transformations.  \emph{Israel J. Math.}  \textbf{26}
(1977), no. 2, 188--196.

\bibitem{KeyNew}
H. Keynes and D. Newton, A ``minimal", non-uniquely ergodic interval exchange transformation.  \emph{Math. Z.}
\textbf{148}  (1976), no. 2, 101--105.


\bibitem{Kir} E. Kirchberg, The classification of purely infinite C*-algebras using Kasparov's theory, to appear in the Fields Institute
Communication series.


\bibitem{KriDim}    W. Krieger, On a dimension for a class of homeo\-morphism groups.  \emph{Math. Ann.}  \textbf{252}
(1979/80), no. 2, 87--95.

\bibitem{KriShift}
W. Krieger,  On dimension functions and topological Markov chains.  \emph{Invent. Math.}  \textbf{56}  (1980), no.
3, 239--250.


\bibitem{KumPrelim} A. Kumjian, Preliminary algebras arising from local homeo\-morphisms.
 \emph{Math. Scand.}  \textbf{52}  (1983),  no. 2, 269--278.

\bibitem{KumFell}  A. Kumjian, Fell bundles over groupoids.
  \emph{Proc. Amer. Math. Soc.}  \textbf{126}  (1998),  no. 4, 1115--1125.


\bibitem{KumPasRae} A. Kumjian, D. Pask, and I. Raeburn, Cuntz-Krieger algebras of directed graphs.
\emph{Pacific J. Math.}  \textbf{184}  (1998),  no. 1, 161--174.

\bibitem{KumRen}
A. Kumjian and J. Renault, KMS states on C*-algebras associated to
expansive maps, \emph{arXiv:math.} OA/0305044 v1 1 May 2003.

\bibitem{MarSevRam} N. Martins, R. Severino, and J. S. Ramos, $K$-theory for Cuntz-Krieger algebras arising from
real quadratic maps. \emph{Int. J. Math. Math. Sci.} \textbf{2003}, no. 34, 2139--2146.

\bibitem{MuhRenWil}
P. S. Muhly, J. N. Renault, and D. P. Williams, Equivalence and isomorphism for groupoid $C\sp
*$-algebras. \emph{J. Operator Theory} \textbf{17} (1987), no. 1, 3--22.

\bibitem{MisSzl}
M. Misiurewicz and W. Szlenk,  Entropy of piecewise monotone mappings. \emph{Studia Math.} \textbf{67} (1980),  no. 1,
45--63.

\bibitem{LinMar}
D. Lind and B. Marcus,
 \emph{An introduction to symbolic dynamics and coding.}
 Cambridge University Press, Cambridge, 1995.

\bibitem{Melo-S} 
W. de Melo and S. van Strien,  
 \emph{One-dimensional dynamics.}  
Ergebnisse der Mathematik und ihrer Grenzgebiete (3) [Results in Mathematics and Related Areas (3)], 25. 
 Springer-Verlag, Berlin, 1993. 


\bibitem{New}  M. Newman, \emph{Integral Matrices},  Pure and Applied Mathematics, Vol. 45.
 Academic Press, New York-London, 1972.

\bibitem{Par}
W. Parry, Symbolic dynamics and transformations of the unit interval. \emph{Trans. Amer. Math. Soc.} \textbf{122} 1966
368--378.

\bibitem{ParTun} W. Parry and S. Tuncel, \emph{Classification problems in ergodic theory.} London Mathematical Society Lecture Note
Series, 67. Statistics: Textbooks and Monographs, 41. Cambridge University Press, Cambridge-New York, 1982.

\bibitem{PasEndo}
W. L. Paschke, The crossed product of a $C\sp{*} $-algebra by an  endomorphism. \emph{Proc. Amer. Math. Soc.} \textbf{80}
(1980),  no. 1, 113--118.

\bibitem{PasKth}
W. L. Paschke, $K$-theory for actions of the circle group on $C\sp{*}  $-algebras. \emph{J. Operator Theory}
\textbf{6} (1981),  no. 1, 125--133.

\bibitem{Phi} N. C. Phillips,  A classification theorem for nuclear purely infinite simple $C\sp *$-algebras.  \emph{Doc.
Math.} \textbf{5}  (2000), 49--114 (electronic).

\bibitem{Pim} M. Pimsner,  A class of C*-algebras generalizing both Cuntz-Krieger algebras and crossed products by $\Z$.
\emph{Free probability theory (Waterloo, ON, 1995)},  189--212,
 Fields Inst. Commun., 12,
 Amer. Math. Soc., Providence, RI, 1997.


\bibitem{PimVoi} M. Pimsner and D. Voiculescu,
 $K$-groups of reduced crossed products by free groups.
 \emph{J. Operator Theory} \textbf{8}  (1982), no. 1, 131--156.

\bibitem{Poo} Y. T. Poon, A $K$-theoretic invariant for dynamical systems, \emph{Trans. Amer. Math. Soc.} \textbf{311}
(1989), 515--533.

\bibitem{Pre}
C. Preston, \emph{Iterates of piecewise monotone mappings on an interval.} Lecture Notes in Mathematics, 1347.
Springer-Verlag, Berlin, 1988.

\bibitem{PutExc1}
I. Putnam, The C*-algebras associated with minimal homeo\-morphisms of the Cantor set, \emph{Pacific J. Math.}
\textbf{136} (1989), no. 2,  329--254.

\bibitem{PutExc2}
I. Putnam,
 C*-algebras arising from interval exchange transformations.
 \emph{J. Operator Theory}  \textbf{27}  (1992),  no. 2, 231--250.


\bibitem{PutSchSka} I. Putnam, K. Schmidt, C. Skau,
 C*-algebras associated with Denjoy homeo\-morphisms of the circle.
 \emph{J. Operator Theory}  \textbf{16}  (1986),  no. 1, 99--126.



\bibitem{RenBook}
J. Renault, \emph{A groupoid approach to $C\sp{*} $-algebras.} Lecture Notes in Mathematics, 793. Springer,
Berlin, 1980.

\bibitem{RenCuntzAlg}
J. Renault, Cuntz-like Algebras, \emph{arXiv:math.} OA/9905185 29 May 1999.

\bibitem{RenRN} J. Renault, The Radon-Nikodym problem for approximately proper equivalence relations,
\emph{arXiv:math.} OA/0261369 v1 23 Nov 2002.

\bibitem{RorBook}
M. R\o rdam, F. Larsen, and N. J. Lausten, \emph{An Introduction to K-theory for C*-algebras.} London Mathematical Society
Student Texts 49, Cambridge University Press, Cambridge, 2000.

\bibitem{RueBook}
D. Ruelle,  \emph{Dynamical zeta functions for piecewise monotone maps of the  interval.} CRM Monograph Series,
4. American Mathematical Society, Providence, RI, 1994.


\bibitem{paper I}
F. Shultz, Dimension groups for interval maps, \emph{arXiv:math.DS/0405466}.

\bibitem{paper 2}
F. Shultz, Dimension groups for interval maps II: The transitive case, \emph{Ergodic Theory
Dynamical Systems}, to appear. \emph{arXiv:math.DS/0405467}.


\bibitem{Spi} J. Spielberg,  Cuntz-Krieger algebras associated with Fuchsian groups.  \emph{Ergodic
Theory Dynam. Systems}  \textbf{13}  (1993),  no. 3, 581--595.

\bibitem{WalBeta}
P. Walters, Equilibrium states for $\beta$-transformations and related transformations, \emph{Math. Z.} \textbf{159}
(1978) 65-88.

\bibitem{WalBook} P. Walters, \emph{An introduction to ergodic theory.} Graduate Texts in Mathematics, 79. Springer-Verlag, New
York-Berlin, 1982.

\end{thebibliography}
\end{document}